
\input amstex
\documentstyle{amsppt}

\TagsOnRight
\TagsAsText
\refstyle{C} 

\newcount\refCount
\def\newref#1 {\advance\refCount by 1
\expandafter\edef\csname#1\endcsname{\the\refCount}}

\newref ALTK 
\newref BOSCH  
\newref HART   
\newref LANG   
\newref MCMU   
\newref MILNE  
\newref MILN   
\newref MORA   
\newref MORT   
\newref MUMFOG 
\newref REES   
\newref SEGAL  
\newref SESH   
\newref ATAEC  
\newref SILV   

\newcount\eqCount
\define\neweqno#1{\global\advance\eqCount by 1%
\expandafter\xdef\csname xyz#1xyz\endcsname{\the\eqCount}\the\eqCount}
\define\referto#1{\csname xyz#1xyz\endcsname}

\define\Part#1{\par\noindent\leavevmode\rlap{\rom{#1}}
   \hglue1.5\parindent}

\define\Case#1#2{\medbreak{\it Case #1\/}:\quad #2
  \hfil\nobreak\vglue\medskipamount\nobreak\noindent\ignorespaces}

\define\Notation#1#2\par{\smallskip\noindent\hangindent1.8cm
\hbox to1.8cm{\hfil$\displaystyle {#1}$\hfil}\ignorespaces
#2\par}

\redefine\AA{\Bbb A}    \define\CC{\Bbb C} 
\define\FF{\Bbb F} \define\GG{\Bbb G}   
   \define\PP{\Bbb P} 
\define\QQ{\Bbb Q}    
\define\ZZ{\Bbb Z}

 \let\a\alpha \let\b\beta \let\g\gamma \let\d\delta 
    \let\l\lambda \let\m\mu
   \let\p\pi \let\r\rho \let\s\sigma \let\t\tau
 \let\f\phi \let\o\omega

     \let\F\Phi
  \let\O\Omega


\define\gM{{\frak M}}
\define\Mbar{\overline{\M}}
\define\Kbar{\overline K}
\define\Ocal{{\Cal O}}
\define\PS{{\sl PS\/}}

\let\<\langle  \let\>\rangle
\let\setminus\smallsetminus

\define\Fix{{\operatorname{Fix}}}
\define\GL{{\operatorname{GL}}}
\define\Hom{{\operatorname{Hom}}}
\define\M{{\operatorname{M}}}
\define\Mor{{\operatorname{Mor}}}
\define\Per{{\operatorname{Per}}}
\define\PGL{{\operatorname{PGL}}}
\define\Proj{{\operatorname{Proj}\,}}
\define\Rat{{\operatorname{Rat}}}
\define\Res{{\operatorname{Res}}}
\define\SL{{\operatorname{SL}}}
\define\Spec{{\operatorname{Spec}\,}}
\define\Sym{{\operatorname{Sym}}}
\define\univ{{\text{univ}}}

\define\FM{{\underline{\M}}}
\define\FRat{{\underline{\Rat}}}
\define\Sets{{\text{\bf Sets}}}
\define\Schemes{{\text{\bf Sch}}}


\topmatter
\title
The space of rational maps on $\PP^1$
\endtitle
\author
Joseph H. Silverman
\endauthor
\affil
Brown University
\endaffil
\address
Mathematics Department
Box 1917
Brown University
Providence, RI 02912 USA
\endaddress
\email
jhs\@math.brown.edu
\endemail
\thanks
Research partially supported by NSF DMS-9424642.
\endthanks
\date
Tuesday, July 30, 1996
\enddate
\abstract
The set of morphisms $\f:\PP^1\to\PP^1$ of degree $d$ is parametrized
by an affine open subset $\Rat_d$ of $\PP^{2d+1}$. We consider the
action of~$\SL_2$ on $\Rat_d$ induced by the {\it conjugation
action\/} of $\SL_2$ on rational maps; that is, $f\in\SL_2$ acts
on~$\f$ via $\f^f=f^{-1}\circ\f\circ f$.  The quotient space
$\M_d=\Rat_d/\SL_2$ arises very naturally in the study of discrete
dynamical systems on~$\PP^1$. We prove that~$\M_d$ exists as an affine
integral scheme over~$\ZZ$, that $\M_2$ is isomorphic to~$\AA^2_\ZZ$,
and that the natural completion of~$\M_2$ obtained using geometric
invariant theory is isomorphic to~$\PP^2_\ZZ$. These results, which
generalize results of Milnor over~$\CC$, should be useful for studying
the arithmetic properties of dynamical systems.
\endabstract
%
\def\1{1} \def\2{2} \def\3{4} \def\4{3} \def\5{5}
\def\6{6} \def\7{7} \def\8{8} \def\9{9}
\toc
\head \1. Notation and summary of results
\endhead
\head \2. The quotient spaces $\M_d$, $\M_d^s$, and $\M_d^{ss}$
\endhead
\head \4. The functors  $\FRat_d$ and $\FM_d$
\endhead
\head \3. Fixed points, periodic points and multiplier systems
\endhead
\head \5. The space $\M_2$ is isomorphic to $\AA^2$
\endhead
\head \6. The completion $\M_2^s$ of $\M_2$
\endhead
\endtoc
\endtopmatter

\document

\head \S\1. Notation and summary of results
\endhead

A rational map $\f:\PP^1\to\PP^1$ of degree~$d$ over a field~$K$ is
given by a pair of homogeneous polynomials
$$
  \f=[F_a,F_b]=
  [a_0X^d+a_1X^{d-1}Y+\cdots+a_dY^d,b_0X^d+b_1X^{d-1}Y+\cdots+b_dY^d]
$$
of degree~$d$ such that~$F_a$ and~$F_b$ have no common roots
(in~$\PP^1(\Kbar)$). This last condition is equivalent to the
condition that
$$
  \Res(F_a,F_b)\ne0,
$$
where the resultant $\Res(F_a,F_b)$ is a certain
bihomogeneous polynomial in the
coefficients $a_0,a_1,\ldots,a_d,b_0,\ldots,b_d$. We will also
frequently write such maps~$\f$ in non-homogeneous form as
$$
  \f(z)=
    \frac{a_0z^d+a_1z^{d-1}+\cdots+a_{d-1}z+a_d}
    {b_0z^d+b_1z^{d-1}+\cdots+b_{d-1}z+b_d}.
$$
\par
We are interested in studying the space of all rational maps
$\PP^1\to\PP^1$ of degree~$d$. These maps are parametrized by the
coefficients of~$F_a$ and~$F_b$, but notice that these are
homogeneous coordinates, since for any non-zero constant~$c$ we have
$[F_a,F_b]=[cF_a,cF_b]$. Thus the space of rational maps of degree~$d$
is the open subset of~$\PP^{2d+1}$ given by the condition
$\Res(F_a,F_b)\ne0$. Notice that this set is an affine variety, since
it is the complement of a hyperplane.

\definition{Definition}
The {\it space of rational maps of degree~$d$} is the affine
open subscheme of $\PP^{2d+1}_\ZZ=\Proj\ZZ[a_0,\ldots,b_d]$ defined by
$$
  \Rat_d = \PP^{2d+1}_\ZZ\setminus\{\Res(F_a,F_b)=0\}.
$$
\enddefinition

To ease notation, we will write
$$\align
  A_d &=\ZZ[a_0,a_1,\ldots,a_d,b_0,b_1,\ldots,b_d],\\
  \rho&=\rho(a,b)=\Res(F_a,F_b)\in A_d.\\
  \endalign
$$
Then $\Rat_d=\Proj A_d\setminus\{\rho=0\}$, so
$$\align
  H^1(\Rat_d,\Ocal_{\Rat_d}) 
  &= A_d[\rho^{-1}]_{(0)} \\
  &= \ZZ\left[\frac{a_0^{i_0}a_1^{i_1}\cdots a_d^{i_d}
      b_0^{j_0}b_1^{j_1}\cdots b_d^{j_d}}{\r}
    \right]_{i_0+\cdots+i_d+j_0+\cdots+j_d=2d},\\
  \endalign
$$
where the ``$(0)$'' subscript denotes elements of degree~$0$ (i.e.,
rational functions whose numerator and denominator are homogeneous of
the same degree).

\remark{Remark}
The space $\Rat_d(\CC)$ of rational maps over the complex numbers has
been studied in some detail. In particular, Segal~\cite{\SEGAL} has
studied the topology of~$\Rat_d(\CC)$ intrinsically and as a subset
of the space of all continuous maps $\PP^1(\CC)\to\PP^1(\CC)$ of
degree~$d$. For example, he proves that the fundamental group
$\p_1\bigl(\Rat_d(\CC)\bigr)$ is cyclic of order~$2d$ and he
gives an explicit description of the universal cover
of~$\Rat_d(\CC)$. We will not consider topological questions of this
nature in this paper.
\endremark

The general linear group $\GL_2$ acts on~$\PP^1$ via linear
fractional transformations in the usual way,
$$
  \pmatrix \a&\b\\ \g&\d\\ \endpmatrix: [X,Y]\longmapsto
  [\a X+\b Y,\g X+\d Y].
$$
The scalar matrices $\left({\a\atop0}\,{0\atop\a}\right)$ act
trivially, so~$\GL_2$ actually acts through its quotient
$\PGL_2=\GL_2/\GG_m$. For various reasons, we will instead consider
the action of the special linear group $\SL_2$. There is very little
lost in doing this, since over an algebraically closed field, the map
$\SL_2\to\PGL_2$ is surjective with kernel equal to $\{\pm1\}$. (In
general over a field, one has
$$
  1@>>>\m_n(K)@>>>\SL_n(K)@>>>\PGL_n(K)@>\det>>K^*/{K^*}^n@>>>1.)
$$
\par
The action of~$\SL_2$ on~$\PP^1$ induces several actions on the
space of rational functions. The one we will be interested in is the
conjugation action given as follows: 
$$\multline
  \text{For $f=\pmatrix \a&\b\\ \g&\d\\ \endpmatrix\in\SL_2$
    and $\f=[F_a,F_b]\in\Rat_d$,}\\
  \f^f = f^{-1}\circ\f\circ f
  = \bigl[ \d F_a(\a X+\b Y,\g X+\d Y)-\b F_b(\a X+\b Y,\g X+\d Y), \\
    -\g F_a(\a X+\b Y,\g X+\d Y)+\a F_b(\a X+\b Y,\g X+\d Y) \bigr].
  \endmultline
$$

\definition{Definition}
The {\it space of conjugacy classes of rational maps of degree~$d$}
is the quotient space of~$\Rat_d$ by the conjugacy action of~$\SL_2$
(in whatever sense this quotient exists). It is denoted by
$$
  \M_d = \Rat_d/\SL_2.
$$
The natural projection map from $\Rat_d$ to $\M_d$ will be denoted
$$
  \<\;\cdot\;\>:\Rat_d\longrightarrow\M_d.
$$
\enddefinition

Our principal aim in this paper is to study the extent to
which~$\M_d$ has any sort of nice structure. A priori, about the
only thing one can say is that for an algebraically closed
field~$\O$, the quotient $\M_d(\O)=\Rat_d(\O)/\SL_2(\O)$ exists as a
set. 

\remark{Remark}
Over the complex numbers, it seems to be known (but not
published?) that $\M_d(\CC)$ has a natural structure as a complex
orbifold, and this is made explict for~$\M_2(\CC)$ in~\cite{\MILN}. 
In fact, Milnor shows that $\M_2(\CC)\cong\CC^2$, and this is one of
the results we will generalize in this paper. See also~\cite{\REES} for
a detailed analysis of various parameter spaces for rational maps of
degree two over~$\CC$ and the loci corresponding to rational maps
which have various complex dynamical properties. 
\endremark

Our first result says that the quotient space~$\M_d$ exists as a
geometric quotient scheme over~$\ZZ$ in the sense of Mumford's
geometric invariant theory. We can further use geometric invariant
theory to deduce various properties about~$\M_d$ and to construct a
natural completion. 

\proclaim{Theorem 1.1}
The quotient $\M_d=\Rat_d/\SL_2$ exists as a geometric quotient
scheme over~$\Spec\ZZ$. It is an affine integral connected scheme
whose affine coordinate ring is the ring of invariant functions
$$
  H^0(\M_d,\Ocal_{\M_d}) = 
  H^1(\Rat_d,\Ocal_{\Rat_d})^{\SL_2} = 
  \left(A_d[\rho^{-1}]_{(0)}\right)^{\SL_2}.
$$
\endproclaim

\remark{Remark}
The precise definition of geometric quotient can be found
in~\cite{\MUMFOG, definition~0.6}. Briefly, in addition to those
properties described in the theorem, the quotient scheme
$\M_d/\ZZ$ has the following pleasant properties:
\roster
\item
The following diagram commmutes:
$$\CD
  \SL_2\times_\ZZ\Rat_d
    @>\text{action of}>\text{$\SL_2$ on $\Rat_d$}>
    \Rat_d \\
  @VV\operatorname{proj}_2 V   @VV\<\,\cdot\,\>V \\
  \Rat_d @>{\<\,\cdot\,\>}>>  \M_d \\
  \endCD
$$
Intuitively, this says that the action of~$\SL_2$ on~$\Rat_d$
descends to the trivial action on~$\M_d$.
\item
For any algebraically closed field~$\O$, the natural map
$\<\,\cdot\,\>:\Rat_d(\O)\to\M_d(\O)$ is surjective, and its fibers
are the $\SL_2(\O)$ orbits of points in $\Rat_d(\O)$.
\item
If $U\subset\M_d$ is an open set, then its inverse image in~$\Rat_d$
is also open.
\endroster
\endremark

As remarked above, Milnor~\cite{\MILN} proved that
$\M_2(\CC)\cong\CC^2$. More precisely, he
describes explicitly two functions $\s_1,\s_2$
on~$\Rat_2=\PP^5\setminus\{\r=0\}$ which are invariant under
the action of~$\SL_2(\CC)$ and which induce a bijection
$(\s_1,\s_2):\M_2(\CC)\to\CC^2$.  We will prove the following
generalization of Milnor's result.

\proclaim{Theorem 1.2}
There are functions
$$
  \s_1,\s_2\in  \left(A_2[\rho^{-1}]_{(0)}\right)^{\SL_2}
$$
(given explicitly in section~\5) which are invariant under the
action of~$\SL_2$ and which induce an isomorphism
$$
  (\s_1,\s_2):\M_2 @>\sim>>\AA^2_\ZZ
$$
of schemes over~$\ZZ$.
\endproclaim

Geometric invariant theory also provides the means to embed~$\M_d$ in
larger quotient spaces, as described in the following result.

\proclaim{Theorem 1.3}
There are open subschemes of~$\PP^{2d+1}$ (over~$\ZZ$)
$$
  \Rat_d \subset (\PP^{2d+1})^s \subset (\PP^{2d+1})^{ss}
$$
which are invariant under the conjugation action
of~$\SL_2$ and such that the quotients
$$
  \M_d=\Rat_d/\SL_2,\qquad \M_d^s=(\PP^{2d+1})^s/\SL_2,
  \qquad\text{and}\qquad \M_d^{ss}=(\PP^{2d+1})^{ss}/\SL_2
$$
exist. More precisely,~$\M_d^s$ is a geometric quotient, $\M_d^{ss}$
is a categorical quotient which is proper and of finite type 
over~$\ZZ$, and~$\M_d$ sits as a dense open subset of both~$\M_d^s$
and~$\M_d^{ss}$.
\endproclaim

The spaces $\M_d^s$ and~$\M_d^{ss}$ are called the spaces of stable
and semi-stable conjugacy classes of rational maps respectively.
Intuitively, the stable locus~$(\PP^{2d+1})^s$ is the largest
set for which the quotient by~$\SL_2$ satisfies
$$
  (\PP^{2d+1})^s(\O)/\SL_2(\O) @>\sim>> \M_d^{s}(\O)
  \qquad\text{for all algebraically closed fields $\O$.}
$$
The semi-stable quotient~$\M_d^{ss}$ has the less agreeable property
that two points in~$(\PP^{2d+1})^{ss}(\O)$ map to the same point
in~$\M_d^{ss}(\O)$ if there is a common point in the closure of
their~$\SL_2(\O)$ orbits in~$\PP^{2d+1}(\O)$. Of course, this weaker
quotient property is balanced by the fact that~$\M_d^{ss}$ is proper
(intuitively, has compact fibers) over~$\Spec\ZZ$. 
\par
In section~\2 we will use Mumford's numerical criterion to describe
(in some sense) the stable and semi-stable loci in~$\PP^{2d+1}$. As a
consequence of that description, we will be able to prove the
following useful result.

\proclaim{Corollary 1.4}
The stable and semi-stable loci coincide if and only if~$d$ is even.
Hence if~$d$ is even, then $\M_d^s=\M_d^{ss}$ is both a geometric
quotient and is proper over $\Spec\ZZ$.
\endproclaim

Working over~$\CC$, Milnor~\cite{\MILN} shows that the
space~$\M_2(\CC)\cong\CC^2$ has a natural compactification
$\hat\M_2(\CC)\cong\PP^2(\CC)$. As Milnor says, this compactification
is natural in the sense that the extra points at infinity ``can be
thought of very roughly as the limits of quadratic rational maps as
they degenerate towards a fractional linear or constant map. However,
caution is needed, since such a limit cannot be uniform over'' all
of~$\PP^1(\CC)$. From the viewpoint of geometric invariant theory,
the ``compactification'' $\M_2^s$ of~$\M_2$ naturally consists
of~$\M_2$, an extra affine line~$\AA^1$, and an extra
point. 

\proclaim{Theorem 1.5}
There is a natural isomorphism $\M_2^s\cong\PP^2$ over~$\ZZ$ so that
the following diagram commutes:
$$\CD
  \M_2 @>\sim>(\s_1,\s_2)> \AA^2 \\
  @VVV   @VVV \\
   \M_2^s @>\sim>> \PP^2\rlap{.} \\
  \endCD
$$
\endproclaim

The functions~$\s_1,\s_2$ are defined in terms of the multipliers
associated to a rational map. We will postpone a complete
definition until section~\3 and be content here to describe
them geometrically. Let~$\O$ be an algebraically closed field, and let
$\f\in\Rat_d(\O)$ be a rational map of degree~$d$ defined over~$\O$.
The fixed points of~$\f$ are the points
$$
  \Fix(\f) = \bigl\{P\in\PP^1(\O)\,:\,\f(P)=P\bigr\}.
$$
We consider this to be a set with multiplicities. Counted with
multiplicty, the set~$\Fix(\f)$ contains exactly $d+1$ points.
If~$P\in\Fix(\f)$, then the derivative
$\f'(P)\in\O$ is well-defined independent of the choice of coordinates
on~$\PP^1$; that is, it depends only on the conjugacy class
$\<\f\>\in\M_d(\O)$. The number~$\f'(P)$ is called the multiplier
of~$\f$ at~$P$. A basic identity asserts that
$$
  \sum_{P\in\Fix(\f)} \frac{1}{1-\f'(P)} = 1.
$$
(See \cite{\MILN} for an analytic proof. But this formula is
essentially algebraic in nature, so the analytic proof implies that
it is a formal identity, hence valid over any field.)
\par
The  individual multipliers form an unordered set, so we
take the corresponding elementary symmetric functions:
$$
  \prod_{P\in\Fix(\f)} \bigl(T+\f'(P)\bigr)
  =\sum_{i=0}^{d+1} \s_i(\f)T^{d+1-i}.
$$
The~$\s_i$'s depend only on the conjugacy class~$\<\f\>$, and their
definition is clearly algebraic, so they give functions on~$\M_d$.
More generally, we can use points of period~$n$,
$$
  \Per_n(\f) = \Fix(\f^n),
$$
and compute the multipliers and symmetric functions of~$\f^n$ at the
points in~$\Per_n(\f)$. These, too, will give functions on~$\M_d$,
which we will denote by~$\s_i^{(n)}$, $i=1,2,\ldots$. (Actually,
it is more efficient to define these functions using
only orbits of formal period~$n$. See section~\3 for the precise
definition of  the~$\s_i^{(n)}$'s.) 
\par
In \cite{\MILN}, Milnor uses his description $\M_2(\CC)\cong\CC^2$ to
show that for maps of degree two, every~$\s_i^{(n)}$ is a polynomial
in $\CC[\s_1,\s_2]$. We can use the above Theorem to strengthen this.

\proclaim{Corollary 1.6}
Every invariant function on~$\M_2$, including in particular
the $\s_i^{(n)}$'s, is a polynomial in $\ZZ[\s_1,\s_2]$.
\endproclaim


\head \S\2. The quotient spaces $\M_d$, $\M_d^s$, and $\M_d^{ss}$
\endhead

In this section we will use geometric invariant theory to construct
the quotient spaces  $\M_d$, $\M_d^s$, and $\M_d^{ss}$. We will
follow closely the methods described in~\cite{\MUMFOG}. We will try
to give complete references to the required results
from~\cite{\MUMFOG}, but in the interest of brevity, we will not take
the time to repeat all of the requisite definitions.
\par
The main construction in~\cite{\MUMFOG} says that if a reductive
group~$G$ acts linearly on a variety (or scheme)~$X$, then the
stable locus $X^s\subset X$ admits a geometric quotient~$Y^s=X^s/G$,
and the semi-stable locus $X^{ss}$ admits a categorical
quotient~$Y^{ss}=X^{ss}/G$. Further, the semi-stable
quotient~$Y^{ss}$ will be proper  (complete) over the base in most
situations. In addition, various nice properties of~$X$ descend to the
quotients~$Y^s$ and~$Y^{ss}$. Applying this general theory to our
specific situation yields the following result.

\proclaim{Theorem \2.1}
We use the notation from  section~\1.
\Part{(a)}
The space of rational function $\Rat_d\subset\PP^{2d+1}$ is an
$\SL_2$-invariant dense open subset of the stable
locus~$(\PP^{2d+1})^s$ in~$\PP^{2d+1}$.  Hence the geometric quotient
$\M_d=\Rat_d/\SL_2$ exists as a scheme over~$\ZZ$.
\Part{(b)}
The geometric quotient
$\M_d^s=(\PP^{2d+1})^s/\SL_2$ and the categorical quotient
$\M_d^{ss}=(\PP^{2d+1})^{ss}/\SL_2$ exist as schemes over~$\ZZ$, and
the natural inclusions
$$
  \M_d\subset \M_d^s\subset \M_d^{ss}
$$
exhibit  each scheme as a dense open subscheme of the next.
\Part{(c)}
The schemes~$\M_d$,~$\M_d^s$, and~$\M_d^{ss}$ are all connected,
integral, normal, and of finite type over~$\ZZ$. Further,~$\M_d$ is
affine and~$\M_d^{ss}$ is proper over~$\ZZ$.
\Part{(d)}
More precisely, if we let $A_d=\ZZ[a_0,\ldots,a_d,b_0,\ldots,b_d]$ and
$\r=\Res(F_a,F_b)\in A_d$, then
$$
  \M_d^{ss}\cong \Proj A_d^{\SL_2} 
  \qquad\text{and}\qquad
  \M_d\cong\Spec A_d[\r^{-1}]_{(0)}^{\SL_2} .
$$
The indicated rings of invariants $A_d^{\SL_2}$ and
$A_d[\r^{-1}]_{(0)}^{\SL_2}$ are finitely generated over~$\ZZ$.
\endproclaim
\demo{Proof}
(a)
The fact that $\Rat_d\subset(\PP^{2d+1})^s$ can be proven similarly
to the proof of \cite{\MUMFOG, proposition~4.2}, using the resultant
form $\r(a,b)=\Res(F_a,F_b)$ in place of the discriminant form.
Alternatively, the inclusion $\Rat_d\subset(\PP^{2d+1})^s$ follows
immediately from the numerical criterion (Proposition~\2.2) proven
below. It is also clear that~$\Rat_d$ is an $\SL_2$-invariant subset
of~$\PP^{2d+1}$, since~$\SL_2$ fixes the resultant form.
Hence~$\Rat_d$ is an
$\SL_2$-stable and $\SL_2$-invariant scheme, so its geometric
quotient exists. Over a field, this is a consquence of Mumford's
construction of quotients \cite{\MUMFOG, chapter~1}, and over~$\ZZ$
it follows by essentially the same methods using  Seshadri's
theorem that a reductive group scheme is geometrically reductive. See
\cite{\SESH} and \cite{\MUMFOG, appendix~1.G}.
\Part{(b)}
The existence of the quotients follows from the work of Mumford and
Seshadri as cited in~(a). The fact that the inclusions are dense open
immersions follows from the analogous fact for the inclusions
$\Rat_d\subset(\PP^{2d+1})^s\subset(\PP^{2d+1})^{ss}$.
\Part{(c,d)}
The schemes $\Rat_d$, $(\PP^{2d+1})^s$, and $(\PP^{2d+1})^{ss}$ are
open subschemes of~$\PP^{2d+1}$, so they are all connected, integral,
and normal. It follows from \cite{\MUMFOG, section~2, remark~(2)}
that the quotients~$\M_d$,~$\M_d^s$, and~$\M_d^{ss}$ have the same
properties. The fact that~$\M_d$ is affine and~$\M_d^{ss}$ is proper
and of finite type over~$\ZZ$ also follows from Seshadri's
work~\cite{\SESH} (see also \cite{\MUMFOG, theorem~1.1 and
appendix~1.G}), as does the description of~$\M_d$ and~$\M_d^{ss}$ via
rings of invariants in~(d).
\enddemo

Next we use Mumford's numerical criterion to describe exactly which
points in~$\PP^{2d+1}$ are (semi)-stable for the action of~$\SL_2$.

\proclaim{Proposition \2.2}
Identifying $\PP^{2d+1}$ with pairs of homogeneous polynomials
$$
  \f=[F_a,F_b]=
  [a_0X^d+a_1X^{d-1}Y+\cdots+a_dY^d,b_0X^d+b_1X^{d-1}Y+\cdots+b_dY^d],
$$
we let~$\SL_2$ act via conjugation as described in section~1. Also
let~$\O$ be an algebraically closed field.
\Part{(a)}
A point in $\PP^{2d+1}(\O)$ is unstable if and only if, after
an $\SL_2(\O)$-conjugation, it
satisfies
$$
  \text{$a_i=0$ for all $i\le\dfrac{d-1}{2}$\quad and\quad
        $b_i=0$ for all $i\le\dfrac{d+1}{2}$.}
$$
\Part{(b)}
A point in $\PP^{2d+1}(\O)$ is not stable if and only if,
after an $\SL_2(\O)$-conjugation, it satisfies
$$
  \text{$a_i=0$ for all $i<\dfrac{d-1}{2}$\quad and\quad
        $b_i=0$ for all $i<\dfrac{d+1}{2}$.}
$$
\endproclaim

As a trivial corollary, we obtain the following useful result.

\proclaim{Corollary \2.3}
If $d$ is even, then every semi-stable point is stable, so
$\M_d^s=\M_d^{ss}$.
\endproclaim

\remark{Remark}
Let $K$ be a non-algebraically closed field of characteristic~0. It
is an interesting question to ask whether the natural map
$$
  \Rat_d(K)\longrightarrow\M_d(K)
$$
is surjective. This is equivalent to asking whether the field of
moduli of a conjugacy class of maps~$\<\f\>$ is also a field of
definition. This question was studied in~\cite{\SILV}, where it is
proved that if~$d$ is even, then $\Rat_d(K)$ always surjects
onto~$\M_d(K)$; but if~$d$ is odd and the Brauer group of~$K$ is
non-trivial, then it never surjects. It is tempting to speculate that
the even/odd dichotomies in~\cite{\SILV} and corollary~\2.3 are
related to one another.
\endremark

\demo{Proof of Proposition \2.2}
We will use the numerical criterion described in \cite{\MUMFOG,
chapter~2}. For a similar computation, see \cite{\MUMFOG, chapter~4,
sections~1 and~2}.
\par
Fix a maximal torus $T\subset\SL_2$. After a change of coordinates,
the action of~$T$ on its canonical representation space~$\AA^2$ can
be diagonalized, so~$T$ becomes the group of matrices
$$
  \pmatrix \a&0\\ 0&\d\\ \endpmatrix
  \qquad\text{subject to the condition $\a\d=1$.}
$$
We are identifying the space of pairs~$[F_a,F_b]$ with the
projective space~$\PP^{2d+1}$, and the canonical action of~$\SL_2$
on~$\AA^2$ is dual to the action on forms, so the (conjugation) action
of an element of $f=\left({\a\atop0}\,{0\atop\d}\right)\in T$
on a point $\f=[F_a,F_b]\in\PP^{2d+1}$ is given by
$$
  \f^f 
  = [F_a,F_b]^f
  = [\a F_a(\a^{-1}X,\d^{-1}Y), \d F_b(\a^{-1}X,\d^{-1}Y) ].
$$
So if we write $F_a=\sum a_iX^{d-i}Y^i$ and $F_b=\sum b_iX^{d-i}Y^i$,
then the action on the $(a,b)$-coordinates is given explicitly as
$$
  a_i\longmapsto \a^{i+1-d}\d^{-i}a_i
  \qquad\text{and}\qquad
  b_i\longmapsto \a^{i-d}\d^{1-i}b_i.
$$
\par
Now consider a one-parameter subgroup (1-\PS) $\l:\GG_m\to\SL_2$.
Attached to each such~$\l$ there is a numerical
invariant~$\m(\f,\l)$ as described in~\cite{\MUMFOG, chapter~2}. The
numerical criterion of~\cite{\MUMFOG, theorem~2.1} says that
$$\align
  \text{$\f$ is unstable}&\Longleftrightarrow
    \text{$\m(\f,\l)<0$ for some 1-\PS\ $\l$,} \\
  \text{$\f$ is not stable}&\Longleftrightarrow
    \text{$\m(\f,\l)\le 0$ for some 1-\PS\ $\l$.} \\
  \endalign
$$
We will now compute this invariant in our situation.
\par
After a change of coordinates, any 1-\PS\ can be
transformed to lie in a maximal torus and be given by
$$
  \l_r(t) = \pmatrix t^r&0\\ 0&t^{-r}\\ \endpmatrix
  \qquad\text{for some integer $r\ge1$.}
$$
The action of~$\l_r$ on~$[F_a,F_b]$ is given by
$$
  a_i\longmapsto t^{-r(d-1-2i)}a_i
  \qquad\text{and}\qquad
  b_i\longmapsto t^{-r(d+1-2i)}b_i.
$$
Then a formula of Mumford~\cite{\MUMFOG, proposition~2.3} says that
$$
  \m(\f,\l_r)=\max\bigl(
    \{r(d-1-2i)\,:\,a_i\ne0\} \cup \{r(d+1-2i)\,:\,b_i\ne0\}\bigr).
$$
Combining this with the numerical criterion says
that~$\f$ is unstable (respectively not stable) if and only if~$\f$ is
conjugate to a map with
$$
  \max\bigl(
    \{r(d-1-2i)\,:\,a_i\ne0\} \cup \{r(d+1-2i)\,:\,b_i\ne0\}\bigr)<0
  \quad\text{(respectively ${}\le0$).}
$$
This is equivalent to the two conditions
$$\gather
  a_i\ne0\Longrightarrow d-1-2i<0 \text{\ (respectively ${}\le0$)}\\
  \qquad\text{and}\qquad\\
  b_i\ne0\Longrightarrow d+1-2i<0 \text{\ (respectively ${}\le0$)},\\
  \endgather
$$
which in turn are the same as
$$\gather
  \text{$a_i=0$ for all $\dfrac{d-1}{2}\ge i$ (respectively ${}>i$)}\\
    \qquad \text{and}\qquad\\
  \text{$b_i=0$ for all $\dfrac{d+1}{2}\ge i$ (respectively ${}>i$).}\\
  \endgather
$$
This completes the proof  proposition~\2.2.
\enddemo

\demo{Proof of Corollary \2.3}
Corollary \2.3 follows immediately from proposition \2.2, since if~$d$
is even, then $(d\pm1)/2$ is not an integer, so the unstable condition
in~(a) and the not-stable condition in~(b) are equivalent.
\enddemo


\head \S\4. The functors  $\FRat_d$ and $\FM_d$
\endhead

In this section we will look at two functors from the category of
schemes to the category of sets. We begin by fixing a realization
$$
  \PP^1=\PP^1_\ZZ= \Proj\ZZ[X,Y].
$$
Equivalently, we fix a basis $X,Y$ for the space of global
sections $H^0(\PP^1,\Ocal_{\PP^1}(1)\bigr)$. 

\definition{Definition}
Let $d\ge1$ be an integer.
The functor $\FRat_d$ of rational maps (really morphisms) of degree~$d$
on~$\PP^1$ is the functor
$$\gather
  \FRat_d:\Schemes\longrightarrow\Sets \\
  \intertext{defined by}
  \FRat_d(S)=\{\text{$S$-morphisms $\f:\PP^1_S\to\PP^1_S$ satisfying
      $\f^*\Ocal_{\PP^1_S}(1)\cong\Ocal_{\PP^1_S}(d)$}\}. \\
  \endgather
$$
\enddefinition

Of course, we still write $\Rat_d$ for the scheme defined in
section~1. That is, $\Rat_d$ is the affine scheme
$$
  \Rat_d=\Spec\ZZ\left[\frac{a_0^{i_0}a_1^{i_1}\cdots a_d^{i_d}
      b_0^{j_0}b_1^{j_1}\cdots b_d^{j_d}}{\r}
    \right]_{i_0+\cdots+i_d+j_0+\cdots+j_d=2d},
$$
where we will write as usual 
$$
  \rho=\rho(a,b)=\Res(F_a,F_b)
$$
for the resultant polynomial. Over $\Rat_d$ we have a universal
morphism $\f^\univ$ of degree~$d$,
$$\matrix
  \PP^1_{\Rat_d} & @>\f^\univ>> & \PP^1_{\Rat_d}. \\
  [X,Y] & @>>>  & [F_a(X,Y),F_b(X,Y)] \\
  \endmatrix
$$

\definition{Definition}
For any scheme~$S$, we define an equivalence relation on the set
$\FRat_d(S)$ as follows. Two $S$-morphisms $\f,\psi\in\FRat_d(S)$ are
equivalent, denoted $\f\sim\psi$, if there is an $S$-isomorphism
$$
  f:\PP^1_S\;\buildrel\sim\over\longrightarrow\PP^1_S
$$
such that $\f\circ f=f\circ\psi$ (i.e., $\f^f=\psi$). We then define
the functor~$\FM_d$ to be the quotient of~$\FRat_d$ by this equivalence
relation:
$$
  \FM_d:\Schemes\longrightarrow\Sets,\qquad
  S\longrightarrow \FRat_d(S)/\sim.
$$
\enddefinition

Our first result says that the functor $\FRat_d$ is representable.

\proclaim{Theorem \4.1}
The scheme $\Rat_d$ represents the functor $\FRat_d$, and in fact the
universal construction described above makes~$\Rat_d$ into a fine
moduli space for~$\FRat_d$.
\endproclaim
\demo{Proof}
Given any $S$-valued point of~$\Rat_d$, say $\s:S\to\Rat_d$, we can
use~$\s$ to base extend the universal map~$\f^\univ$ and obtain a
morphism $\f^\univ_\s\in\FRat_d(S)$ defined by the following diagram:
$$\CD
   \PP^1_{\Rat_d}\times_{\Rat_d} S 
      @>\f^\univ\times\bold1_S>> 
      \PP^1_{\Rat_d}\times_{\Rat_d} S \\
    @| @| \\
  \PP^1_S @>\f^\univ_\s>> \PP^1_S \\
  \endCD
$$
This gives a map from $\Rat_d(S)=\Hom(S,\Rat_d)$ to $\FRat_d(S)$ for
every scheme~$S$.
\par
Next suppose that we start with an element $\f\in\FRat_d(S)$. We will
further suppose that~$S$ is affine, $S=\Spec B$. We have assumed that
we have fixed a $\ZZ$-basis $X,Y$ for
$H^0(\PP^1,\Ocal_{\PP^1}(1)\bigr)$, so~$X,Y$ will certainly be a
$B$-basis for $H^0(\PP_S^1,\Ocal_{\PP_S^1}(1)\bigr)$, and similarly
$X^d,X^{d-1}Y,\ldots,Y^d$ is a (canonical, given the initial choice
of~$X$ and~$Y$) $B$-basis for
$H^0(\PP_S^1,\Ocal_{\PP_S^1}(d)\bigr)$. Further, the definition of
$\FRat_d$ implies that
$\f^*\Ocal_{\PP_S^1}(1)\cong\Ocal_{\PP^1_S}(d)$, so taking global
sections we can write
$$\align
  \f^*X&=\a_0X^d+\a_1X^{d-1}Y+\cdots+\a_dY^d=F_\a,\\
  \f^*Y&=\b_0X^d+\b_1X^{d-1}Y+\cdots+\b_dY^d=F_\b,\\
  \endalign
$$
where $\a_0,\ldots,\b_d\in B$ are uniquely determined by~$\f$ (and
our choice of~$X,Y$).
\par
We are going to prove that $\rho(\a,\b)=\Res(F_\a,F_\b)$ is a unit
in~$B$. Assuming this, we see that $(\a,\b)$ defines a point
$\t=\t(\f)$ in
$\Rat_d(B)$, and then it is clear that $\f^\univ_\t$ is just the
original map~$\f$. This means that the maps
$$\matrix
  \Rat_d(S) & \longrightarrow & \FRat_d(S) 
     &\quad\text{and}\quad&  \FRat_d(S) & \longrightarrow & \Rat_d(S)
       \\
  \s & \longmapsto & \f^\univ_\s 
     && \f & \longmapsto & \t(\f)=(F_{\a(\f)},F_{\b(\f)}) \\
  \endmatrix
$$
are inverses, at least on affine schemes~$S$. However, the uniqueness
of the~$(F_\a,F_\b)$ associated to a given $\f\in\FRat_d(S)$ means
that we can glue to get the same result for arbitrary schemes~$S$.
\par
It remains to show that $\rho(\a,\b)\in B^*$. We know that
$\f:\PP^1_S\to\PP^1_S$ is a morphism, so in particular $\f^*X,\f^*Y$
generate the sheaf $\f^*\Ocal_{\PP^1_S}(1)$. Hence for any
(closed) point $P\in\PP^1_S$, at least one of the sections $(\f^*X)_P$
and $(\f^*Y)_P$ must be non-zero. In other words, for every (maximal)
ideal $\goth p\in\Spec B$, the forms
$$
  \overline{\f^*X}=F_\a(X,Y)\bmod\goth p
  \qquad\text{and}\qquad \overline{\f^*Y}=F_\b(X,Y)\bmod\goth p
$$
define the trivial locus in $\PP^1_{B/\goth p}$. This implies that
$\Res(F_\a,F_\b)\notin\goth p$, and since this is true for all
(maximal) ideals, we conclude as desired that $\Res(F_\a,F_\b)\in
B^*$.
\enddemo

Next we consider the functor~$\FM_d$. 

\proclaim{Theorem \4.2}
There is a natural map of functors
$$
  \FM_d\longrightarrow\Hom(\,\cdot\,,\M_d)
$$
with the property that $\FM_d(\Omega)\cong\M_d(\Omega)$ for every
algebraically closed field~$\Omega$.
\endproclaim
\demo{Proof}
Let $\xi\in\FM_d(S)$. Within the equivalence class~$\xi$ we choose an
element $\f\in\FRat_d(S)$. From theorem~\4.1, we may regard~$\f$ as
an element of $\Rat_d(S)$, and then the construction of~$\M_d$ as a
quotient (theorem~\2.1) gives us a point $\l=\l(\xi,\f)\in\M_d(S)$.
We claim that~$\l$ is independent of the choice of~$\f$, and so gives
a well-defined map $\FM_d(S)\to\M_d(S)$. To verify this, let
$\psi=\f^f$ be another element of~$\xi$, where $f:\PP^1_S\to\PP^1_S$
is an $S$-isomorphism. Then~$\f$ and~$\psi$ are $S$-valued points
of~$\Rat_d$, and we want to show that the compositions
$$
  S  \;\vcenter{\openup-2\jot\halign{\hfil$#$\hfil\cr
    \scriptstyle\phi\cr\longrightarrow\cr
    \longrightarrow\cr\scriptstyle\psi\cr}}\;
  \Rat_d \longrightarrow \M_d
$$
give the same map $S\to\M_d$. Covering~$S$ by affine open sets, we
may assume that $S=\Spec B$.
\par
The $S$-isomorphism~$f$ satisfies
$f^*\Ocal_{\PP_S^1}(1)\cong\Ocal_{\PP_S^1}(1)$, so
$$
  f^*X=\a X+\b Y\qquad\text{and}\qquad
  f^*Y=\g X+\d Y
$$
for some $\a,\b,\g,\d\in B$. Further, the fact that~$f$ has an
inverse means that $\det f=\a\d-\b\g\in B^*$. Let
$B'=B\bigl[\sqrt{\a\d-\b\g}\bigr]$ and $S'=\Spec B'$. Notice that
$B'/B$ is a finite extension, so
$S'\to S$ is surjective. This allows us to replace~$S$ by~$S'$, and
then we may replace~$f$ with the map~$f'$ determined by the conditions
$$\define\[#1]{\frac{#1}{\sqrt{\a\d-\b\g}}}
  {f'}^*X=\[\a] X+\[\b] Y\quad\text{and}\quad
  {f'}^*Y=\[\g] X+\[\d] Y.
$$
It is still true that $\phi\circ f'=f'\circ\psi$, and now $\det f'=1$.
\par
Thus $\f$ and $\psi$ are $\SL_2(B')$-equivalent, so any function in the
ring of invariants
$$
  H^0(\M_d,\Ocal_{\M_d}) = H^0(\Rat_d,\Ocal_{\Rat_d})^{\SL_2}
$$
will take the same value at~$\f$ and~$\psi$. Since~$\M_d$ is the
spectrum of this ring, it follows that~$\f$ and~$\psi$ give the same
$S$-valued point of~$\M_d$. This completes the proof that the map
$\FM_d(S)\to\M_d(S)$ defined above is indeed well defined,
independent of the choice of a representative in $\Rat_d(S)$.
\par
We also need to show that $\FM_d(S)$ is isomorphic to $\M_d(S)$ on
geometric points $S=\Spec\Omega$ (i.e., where $\Omega$ is an
algebraically closed field). But this is clear, since over an
algebraically closed field we have
$$
  \FM_d(\Omega) = \Rat_d(\Omega)/\PGL_2(\Omega)
  \qquad\text{and}\qquad
  \M_d(\Omega) = \Rat_d(\Omega)/\SL_2(\Omega),
$$
and the map $\SL_2(\Omega)\to\PGL_2(\Omega)$ is surjective, so the
quotients are the same.
\enddemo


\head \3. Fixed points, periodic points and multiplier systems
\endhead
In this section we will construct functions on the quotient space
$\M_d$ by associating to each $\f\in\M_d$ the system of multipliers
of its fixed (or more generally periodic) points. To motivate this
construction, we begin by describing the situation over an
algebraically closed field~$k$. 
\par
Thus let $\f\in\Rat_d(k)$, so $\f:\PP^1_k\to\PP^1_k$
is a rational map of degree~$d$. Such a
map has exactly $d+1$ fixed points (counted with multiplicity), say
$\xi_1,\ldots,\xi_{d+1}$. For each $\xi=\xi_i$, the map $\f$ induces a 
$k$-linear map~$\f^*$ from~$\O_{\PP_k^1,\xi}$ to itself, where
$\O_{\PP_k^1,\xi}$ denotes the space of germs of differential 1-forms
at~$\xi$. This vector space has dimension~1 over~$k$, so if
we take any non-zero differential form~$\o\in\O_{\PP_k^1,\xi}$, then
$\f^*(\o)$ is a multiple of~$\o$, say
$\f^*(\o)=\f'(\xi)\o$. The number $\f'(\xi)\in k$ is called the {\it
multiplier of~$\f$ at the fixed point~$\xi$}, and the set
$\{\f'(\xi_1),\f'(\xi_2),\ldots,\f'(\xi_{d+1})\}$ is the associated
{\it multiplier system}. 
\par
The fixed points depend algebraically, but not rationally, on the
coefficients of~$\f$,  and in any case they only form an unordered
set. We define quantities $\s_i=\s_i(\f)$ by taking the symmetric
functions of the multipliers:
$$
  \prod_{i=1}^{d+1}\bigl(T+\f'(\xi_i)\bigr)
  =\sum_{i=0}^{d+1} \s_i T^{d+1-i}.
$$
The $\s_i$'s are symmetric functions of the $\f'(\xi_i)$'s, and each
$\f'(\xi)$ is a rational function of~$\xi$, so the $\s_i$'s are
actually rational functions of the coefficients of~$\f$. In other
words, the $\s_i$'s are rational (in fact, regular) functions
on~$\Rat_d/k$.
\par
Now let $f\in\PGL_2(k)$ be an automorphism of $\PP^1_k$. Then
$f^{-1}(\xi_1),\ldots,f^{-1}(\xi_{d+1})$ are the fixed points of 
$\f^f$, and the chain rule tells us that
$$\align
  \left(\f^f\right)'\bigl(f^{-1}(\xi)\bigr)
  &=\left(f^{-1}\circ\f\circ f\right)'\bigl(f^{-1}(\xi)\bigr) \\
  &=(f^{-1})'(\f(\xi))\cdot\f'(\xi)\cdot f'(f^{-1}(\xi)) \\
  &=(f^{-1})'(\xi)\cdot\f'(\xi)\cdot f'(f^{-1}(\xi)) \\
  &=\f'(\xi). \tag\neweqno{3.2} \\
  \endalign
$$
Thus the multiplier system of~$\f$ is $\PGL_2(k)$-invariant, so the
$\s_i$'s descend to give regular functions on~$\M_d$.
\par
We now explain how this construction generalizes over~$\ZZ$.

\proclaim{Theorem \3.1}
Let $\f=\f^\univ:\PP^1_{\Rat_d}\to\PP^1_{\Rat_d}$ be the universal
morphism of degree~$d$ as described in Section~\4. There exists a
unique reduced closed subscheme
$$
  \Fix\subset\PP^1_{\Rat_d}
$$
having the following two properties:
\Part{(i)}
$\f\big|_{\Fix}={\bold1}_{\Fix}$, i.e., $\f$ induces the identity
map on~$\Fix$.
\Part{(ii)}
If $Z\subset\PP^1_{\Rat_d}$ is a reduced closed subscheme with the
property that $\f\big|_Z={\bold1}_Z$, then $Z\subset\Fix$.
\par
The subscheme $\Fix$ also satisfies:
\Part{(iii)}
$\Fix$ is integral (i.e., reduced and irreducible).
\Part{(iv)}
The projection $\Fix\to\Rat_d$ is a finite morphism of degree $d+1$.
\endproclaim
\demo{Proof}
To ease notation, we will write $\f=\f^\univ=[F_a,F_b]$. If $\Fix$
exists, its uniqueness is clear from~(i) and~(ii), so we just need to
find a subscheme with properties~(i) and~(ii). We set
$$
  \Fix = V\bigl(YF_a(X,Y)-XF_b(X,Y)\bigr)\subset\PP^1_{\Rat_d}.
  \tag\neweqno{3.1}
$$
In other words, we start with the hypersurface of type~$(d+1,1)$
in $\PP^1\times\PP^{2d+1}$ defined by the bihomogeneous form
$YF_a-XF_b$, and then we take its intersection with
$\PP^1\times\Rat_d$. It is clear that~$\f$ fixes~$\Fix$,
since~$\Fix$ is the subscheme defined by the ``condition''
$\f([X,Y])=[X,Y]$. 
\par
Next let $Z\subset\PP^1_{\Rat_d}$ be fixed by~$\f$. If $\s:\Spec k\to
Z$ is any geometric point of~$Z$, say $\s(k)=[\a,\b]\in\PP^1(k)$,
then $\f$ fixes $\s(k)$, so
$$
  [\a,\b]=\f([\a,\b])=[F_a(\a,\b),F_b(\a,\b)].
$$
Hence $\s(k)\in\Fix(k)$, which shows that every geometric point of~$Z$
lies in~$\Fix$. It follows  that $Z$ is a subscheme 
of~$\Fix$ (this is where we use the assumption that~$Z$ is reduced).
This completes the proof that the subscheme
$\Fix$ defined by~(\referto{3.1}) satisfies~(i) and~(ii).
\par
Next we observe that the bihomogeneous form $YF_a-XF_b$ is clearly
irreducible, since it has degree~1 in~$(a,b)$. More precisely, if it
were to factor in $\ZZ[a,b][X,Y]$, then one of the factors would have
to lie in $\ZZ[X,Y]$, and its clear that $YF_a-XF_b$ has no such
factors. Hence~$\Fix$ is irreducible, and it is reduced by
assumption, which verifies~(iii).
\par
To check~(iv), we observe from~(\referto{3.1}) that $\Fix$ is
clearly quasi-finite over~$\Rat_d$, and similarly it is proper (even
projective) over~$\Rat_d$. Therefore~$\Fix$ is finite over~$\Rat_d$
(see \cite{\HART, exercise~11.2} or \cite{\MILNE, chapter~I,
proposition~1.10}). The degree of the map is then clear, since
$YF_a-XF_b$ is homogeneous of degree $d+1$ in $(X,Y)$.
\enddemo

It follows from Theorem~\3.1(i) that $\f$ induces an
$\Ocal_\Fix$-linear map $\f^*$ from
$$
  \O_{\PP^1_{\Rat_d}/\Rat_d}\otimes_{\Ocal_{\PP^1_{\Rat_d}}}\Ocal_\Fix
$$
to itself. This sheaf is (locally) free of rank~1 over~$\Ocal_\Fix$,
so from~(iii) it is (locally) free of rank $d+1$ over
$\Ocal_{\Rat_d}$. Thus~$\f^*$ defines an~$\Ocal_{\Rat_d}$-linear
map of a (locally) free sheaf of rank $d+1$, so we can
compute its characteristic polynomial
$$
  \det(T+\f^*) = \sum_{i=0}^{d+1}\s_iT^{d+1-i}
$$
to obtain (local) sections $\s_1,\ldots,\s_{d+1}$
of~$\Ocal_{\Rat_d}$.

\proclaim{Proposition \3.2}
The functions $\s_1,\ldots,\s_{d+1}$ described above are global
sections of $\Ocal_{\Rat_d}$. Further, they are invariant under the
conjugation action of~$\SL_2$, and hence they descend to give global
sections of $\Ocal_{\M_d}$.
\endproclaim
\demo{Proof}
The scheme $\Rat_d$ is affine, say $\Rat_d=\Spec B$; and $\Fix$ is
finite over $\Rat_d$, so it too is affine, say $\Fix=\Spec B'$, where
$B'/B$ is a finite extension of degree $d+1$. Then 
$$
  \O_{\PP^1_{\Rat_d}/\Rat_d}\otimes_{\Ocal_{\PP^1_{\Rat_d}}}\Ocal_\Fix
  = \Ocal_{\PP^1_B/B}\otimes_{\Ocal_{\PP^1_B}} B'
$$
is a free $B'$-module of rank~1 generated by~$dz$, where~$z$ is a
uniformizer on~$\PP^1_B$ at~$\Fix$. (We can think of $\Fix$, an
integral subscheme of $\PP^1_B$, as a point of $\PP^1_B$.) The
map~$\f$ induces an endomorphism~$\f^*$ of this module, so 
$$
  \f^*(dz)=\f'\cdot dz\quad\text{for some element $\f'\in B'$.}
$$
Now $B'$ is a free $B$-module of rank $d+1$, so multiplication
by~$\f'$ gives an $B$-linear endomorphism of $B^{d+1}$. The
characteristic polynomial of this endomorphism is well-defined
independent of the choice of a basis, which gives
$$
  \det(T+\f')=\sum_{i=0}^{d+1} \s_iT^{d+1-i}
  \quad\text{for elements $\s_i\in B$.}
$$
In other words, $\s_1,\ldots,\s_{d+1}$ are global sections
of~$\Ocal_{\Rat_d}$.
\par
In order to show that the~$\s_i$'s descend to~$\M_d$, we must show
that they are $\SL_2$-invariant. Since~$\M_d$ is reduced, it suffices
to check invariance on geometric points, so let~$k$ be an
algebraically closed field, let $\f\in\Rat_d(k)$, and let
$f\in\SL_2(k)$. Then the equality $\s_i(\f)=\s_i(\f^f)$ follows from
the chain-rule calculation~(\referto{3.2}). Hence the~$\s_i$'s are
global sections of $(\Ocal_{\Rat_d})^{\SL_2}=\Ocal_{\M_d}$.
\enddemo

We continue to let $\f=\f^\univ$ be the universal morphism of
degree~$d$ over~$\Rat_d$. Theorem~\3.1 above describes the fixed
subscheme of~$\f$. More generally, for any $n\ge1$, we can consider
the periodic subscheme of period~$n$, as described in the following
theorem.

\proclaim{Theorem \3.3}
For every $n\ge1$, there exists a unique reduced closed subscheme
$$
  \Per_n\subset\PP^1_{\Rat_d}
$$
having the following two properties:
\Part{(i)}
$\f^n\big|_{\Per_n}={\bold1}_{\Per_n}$, i.e., $\f^n$ induces the
identity map on~$\Per_n$.
\Part{(ii)}
If $Z\subset\PP^1_{\Rat_d}$ is a reduced closed subscheme with the
property that $\f^n\big|_Z={\bold1}_Z$, then $Z\subset\Per_n$.
\par
The scheme $\Per_n$ is called the {\rm scheme of periodic points of
period~$n$}. The subscheme $\Per_n$ also satisfies:
\Part{(iii)}
The projection $\Per_n\to\Rat_d$ is a finite morphism of degree
$d^n+1$.
\endproclaim
\demo{Proof}
Most of the proof is very similar to the proof of Theorem~\3.1, so we
just briefly sketch. We can write $\f^n=[F_a^{(n)},F_b^{(n)}]$, where
$F_a^{(n)}$ and $F_b^{(n)}$ are bihomogeneous polynomials in
$\ZZ[a,b][X,Y]$ of bidegree $((d^n-1)/(d-1),d^n)$. 
Consider the closed subscheme of~$\PP^1_{\Rat_d}$ defined by the
equation
$$
  YF_a^{(n)}-XF_b^{(n)},
  \tag\neweqno{3.3}
$$
and let $\Per_n$ be this subscheme with the induced reduced subscheme
structure. It is clear that $\f^n$ induces  the identity map
on~(\referto{3.3}), hence also on~$\Per_n$, and then~(i) and~(ii) and
the fact that $\Per_n$ is finite over~$\Rat_d$ are proven in the same
way as Theorem~\3.1. It remains to show that the degree of $\Per_n$ 
over~$\Rat_d$ is exactly $d^n+1$.
\par
Since $\Per_n$ is the subscheme of $\PP^1_{\Rat_d}$ given
by~(\referto{3.3}) with the induced reduced subscheme structure,
and~(\referto{3.3})  has degree $d^n+1$ in the variables~$(X,Y)$, we
must show that the polynomial $YF_a^{(n)}-XF_b^{(n)}$ has no repeated
factors when factored in $\ZZ[a,b][X,Y]$.
For this, it suffices to show that it has no repeated factors when we
specialize~$(a,b)$. Consider the rational map
$\f=[X^d,Y^d]\in\Rat_d(\CC)$. For this map we have
$\f^n=[X^{d^n},Y^{d^n}]$, so~(\referto{3.3}) becomes
$$
  YF_a^{(n)}-XF_b^{(n)} = XY(X^{d^n-1}-Y^{d^n-1}),
$$
which is a polynomial with distinct roots in~$\PP^1(\CC)$. This
proves that~(\referto{3.3}) has no repeated factors in
$\ZZ[a,b][X,Y]$ (and in fact, no repeated factors in
$\FF_p[a,b][X,Y]$ provided $d^n\not\equiv1\pmod{p}$).
\enddemo

Theorem~\3.1 included the assertion that~$\Fix$ is
irreducible, but the analogous statement for~$\Per_n$ was omitted in
Theorem~\3.3. In fact, $\Per_n$ is always reducible for $n\ge2$,
since in particular we always have $\Fix\subset\Per_n$. 
In terms of polynomials, it is easy to check that the
equation~(\referto{3.3}) defining~$\Per_n$  is divisible by
$YF_a-XF_b$. More generally, we can decompose~$\Per_n$ into  pieces as
described in the following theorem.

\proclaim{Thereom \3.4}
With notation as in Theorems~\3.1 and~\3.3, there are unique
reduced closed subschemes
$$
  \Per_m^*\subset\PP^1_{\Rat_d},
$$
one for each $m\ge1$, with the following properties:
\Part{(i)}
$\Per_1^*=\Fix$.
\Part{(ii)}
$\Per_n=\bigcup_{m|n}\Per_m^*$ for every $n\ge1$.
\par
The scheme $\Per_m^*$ is called the {\rm scheme of periodic points of
formal\footnote{In \cite{\MORT}, these were called points
of ``essential'' period $m$ and were denoted $Z_m^*$; but we feel that
Milnor's ``formal'' \cite{\MILN} is a better terminology.} period~$m$}.
In addition:
\Part{(iii)}
$\f^m\big|_{\Per_m^*}={\bold1}_{\Per_m^*}$.
\Part{(iv)}
$\Per_m^*$ is finite over $\Rat_d$, and if we let $\nu_m$ be the
degree of $\Per_m^*$ over $\Rat_d$, then
$$
  d^n+1=\sum_{m|n}\nu_m\qquad\text{and}\qquad
  \nu_m=\sum_{r|m}\mu(m/r)(d^r+1),
$$
where $\m$ is the M\"obius function.
\endproclaim
\demo{Proof}
The proof is by induction on~$m$. We have $\Per_1^*=\Fix$ from~(i),
so we are okay for $m=1$. Now suppose that we know the theorem for
all $m<n$. Consider the scheme $\Per_n$ from Theorem~\3.3. For any
$m|n$ with $m<n$, we know that $\f^m$ fixes $\Per_m^*$, and so
$\f^n=(\f^m)^{(n/m)}$ also fixes $\Per_m^*$. It follows from
Theorem~\3.3 that $\Per_m^*$ is  a subscheme of $\Per_n$. Since
$\Per_n$ and the $\Per_m^*$'s are finite over the (affine
irreducible) scheme~$\Rat_d$, it follows that
$$
  \Per_n = \bigg(\bigcup_{m|n,\,m<n}\Per_m^*\bigg) \cup \Per_n^*,
$$
where $\Per_n^*$ is a union of irreducible components of~$\Per_n$. 
Further, since~$\f^n$ induces the identity map on~$\Per_n$, it
clearly induces the identity map on~$\Per_n^*$; and since $\Per_n$
is finite over $\Rat_d$, the same is true of $\Per_n^*$. Finally, the
degree of $\Per_n^*$ over $\Rat_d$ satisfies
$$\align
  \nu_n&=\deg(\Per_n^*\to\Rat_d)\\
  &=\deg(\Per_n\to\Rat_d)-\sum_{m|n,\,m<n}\deg(\Per_m^*\to\Rat_d)\\
  &=(d^n+1)-\sum_{m|n,\,m<n}\nu_m.\\
  \endalign
$$
This gives the first part of~(iv), and the second part is just
M\"obius inversion.
\enddemo

\remark{Remark}
The scheme $\Per_n\subset\Rat_d$ is given by the vanishing
of the homogeneous polynomial
$$
  \F_n \;\buildrel\text{def}\over=\; YF_a^{(n)}-XF_b^{(n)}.
$$
The defining property of~$\Per_n$ (or a direct calculation) shows that
if $m|n$, then $\F_m|\F_n$ in $\ZZ[a,b][X,Y]$. Then the fact
that~$\F_n$ is reduced (i.e., has no repeated factors) and a simple
inclusion-exclusion argument shows that the product
$$
  \F_n^* \;\buildrel\text{def}\over=\;
  \prod_{m|n} (\F_m)^{\mu(n/m)}.
$$
is in~$\ZZ[a,b][X,Y]$. Looking at the defining properties of~$\Per_n^*$
in Theorem~\3.4, we see that $\Per_n^*$ is given by the equation
$\F_n^*=0$.
\par
It is almost certainly the case that $\F_n^*$, and thus $\Per_n^*$,
are irreducible, but this has not yet been proven. A similar
problem for the space of (monic) polynomial maps is treated by Morton
in~\cite{\MORA}.
\endremark

\remark{Remark}
The scheme $\Per_n^*$ is finite over $\Rat_d$, and its degree $\nu_n$
gives the number of periodic points of formal period~$n$ for a
rational map of degree~$d$. The following table gives the value
of~$\nu_n$ for small values of~$d$ and~$n$.
$$\vcenter{\offinterlineskip
  \def\HRULE{\noalign{\hrule}}
  \halign{\strut#&\vrule#&$\;\hfil#\;$
    &\vrule\hglue1.5pt\vrule#&$\;\hfil#\;$
    &&\vrule#&$\;\hfil#\;$\cr
  \HRULE
  &&d\,\backslash\, n\hfil&&1\hfil&&2\hfil&&3\hfil&&4\hfil
     &&5\hfil&&6\hfil&&7\hfil&&8&\cr
  \HRULE\noalign{\vskip1.5pt}\HRULE
  &&2\hfil&&3&&2&&6&&12&&30&&54&&126&&240&\cr
  \HRULE
  &&3\hfil&&4&&6&&24&&72&&240&&696&&2184&&6480&\cr
  \HRULE
  &&4\hfil&&5&&12&&60&&240&&1020&&4020&&16380&&65280&\cr
  \HRULE
  &&5\hfil&&6&&20&&120&&600&&3120&&15480&&78120&&390000&\cr
  \HRULE
  &&6\hfil&&7&&30&&210&&1260&&7770&&46410&&279930&&1678320&\cr
  \HRULE
  &&7\hfil&&8&&42&&336&&2352&&16800&&117264&&823536&&5762400&\cr
  \HRULE
  &&8\hfil&&9&&56&&504&&4032&&32760&&261576&&2097144&&16773120&\cr
  \HRULE
  \noalign{\vskip2\jot}
  \multispan{19}\hfil The degree of $\Per_n^*$ over $\Rat_d$\hfil\cr
  } }
$$
\endremark

We can define functions using $\Per_n$ and $\Per_n^*$ in exactly the
same way that we defined functions using $\Per_1=\Fix$. Following
Milnor~\cite{\MILN}, we will use the more intrinsic scheme~$\Per_n^*$
of  periodic points of formal period~$n$. Let $\nu_n$ be the degree
of $\Per_n^*$ over $\Rat_d$. Then
$$
  \O_{\PP^1_{\Rat_d}/\Rat_d}\otimes_{\Ocal_{\PP^1_{\Rat_d}}}
  \otimes\Ocal_{\Per_n^*}
$$
is a (locally) free sheaf of rank~1 over $\Ocal_{\Per_n^*}$,
hence (locally) free of rank~$\nu_n$ over $\Ocal_{\Rat_d}$. The
map~$\f$ induces a linear endomorphism of this sheaf, and we compute
the characteristic polynomial
$$
  \det(T+\f^*)=\sum_{i=0}^{\nu_n}\s_i^{(n)}T^{\nu_n-i}
$$
for certain  sections $\s_i^{(n)}$ of $\Ocal_{\Rat_d}$. The following
result generalizes Proposition~\3.2.

\proclaim{Theorem \3.5}
The functions 
$$
  \s_i^{(n)},\qquad \text{for all $n\ge1$ and $1\le i\le\nu_n$,}
$$
as described above, are global sections of
$\Ocal_{\Rat_d}$. Further, they are invariant under the conjugation
action of $\SL_2$, and hence they descend to give global sections
of~$\Ocal_{\M_d}$.
\endproclaim
\demo{Proof}
The proof is the same, mutatis mutandis, as the proof of
Proposition~\3.2.
\enddemo


\head \S\5. The space $\M_2$ is isomorphic to $\AA^2$
\endhead

In this section we will prove the following theorem.

\proclaim{Theorem \5.1}
The natural map
$$
  \M_2\longrightarrow \Spec\ZZ[\s_1,\s_2]\cong\AA^2_\ZZ
$$
is an isomorphism of schemes over~$\ZZ$.
\endproclaim

\remark{Remark}
Theorem~\5.1 may be compared with Milnor's result \cite{\MILN,
lemma~3.1} which says that there is an algebraic bijection between
$\M_2(\CC)$ and~$\CC^2$. Milnor uses his result to deduce \cite{\MILN,
lemma~D.1} that the higher order invariants $\s_i^{(n)}$ are
in~$\CC[\s_1,\s_2]$. He illustrates this corollary with the examples
$$\align
  \s_1^{(2)}&=2\s_1+\s_2,\\
  \s_1^{(3)}&=\s_1(2\s_1+\s_2)+3\s_1+3,\\
  \s_2^{(3)}&=(\s_1+\s_2)^2(2\s_1+\s_2)-\s_1(\s_1+2\s_2)+12\s_1+28.\\
  \endalign
$$
Using Theorem \5.1, we can strengthen Milnor's result by showing
that the $\s_i^{(n)}$'s are always polynomials in~$\s_1$ and~$\s_2$
with {\it integer} coefficients.
\endremark

\proclaim{Corollary \5.2}
The ring of $\SL_2$-invariant functions on $\Rat_2$ is
exactly $\ZZ[\s_1,\s_2]$. In particular, all of the higher order
invariants $\s_i^{(n)}$ are in~$\ZZ[\s_1,\s_2]$.
\endproclaim

\remark{Remark}
If we write
$$
  \f(z)={a_0z^2+a_1z+a_2\over b_0z^2+b_1z+b_2}
  = {F_a(z)\over F_b(z)},
$$
then the corresponding resultant is
$$\align
  \rho&=\rho(a,b)=\Res(F_a,F_b) \\
  &=a_{2}^{2} b_{0}^{2} - a_{1} a_{2} b_{0} b_{1} + 
  a_{0} a_{2} b_{1}^{2} + a_{1}^{2} b_{0} b_{2} - 
  2 a_{0} a_{2} b_{0} b_{2} - a_{0} a_{1} b_{1} b_{2} + 
  a_{0}^{2} b_{2}^{2}. \\
  \endalign
$$
The space of rational functions $\Rat_2$ is the subset of
$\PP^5=\Proj\ZZ[a_0,a_1,a_2,b_0,b_1,b_2]$ given by the  non-vanishing
condition $\rho(a,b)\ne 0$, so $\Rat_2$ is the affine scheme
$$
  \Rat_2 = 
  \Spec  A_2[\r^{-1}]_{(0)}  =
  \Spec \ZZ\left[ 
  {a_0^{i_0}a_1^{i_1}a_2^{i_2}b_0^{j_0}b_1^{j_1}b_2^{j_2}
      \over\rho(a,b)}\right]_{i_0+i_1+i_2+j_0+j_1+j_2=4} 
$$
\par
The action of~$\SL_2$ on~$\Rat_2$ is given by its action on
the~$a_i$'s and~$b_i$'s corresponding to the rule
$\f^f=f^{-1}\circ\f\circ f$. We will omit giving the action
explicitly, but we note from Theorem~\2.1 that~$\M_2$ is the
affine scheme whose affine coordinate ring is the ring of invariants
of this $\SL_2$-action. According to Theorem~\5.1, this ring of
invariants is exactly $\ZZ[\s_1,\s_2]$, so it seems worthwhile to
write down~$\s_1$ and~$\s_2$ explicitly in terms of the~$a_i$'s
and~$b_i$'s. 
$$\align
  \rho(a,b)\s_1(\f) &=a_{1}^{3} b_{0} - 4 a_{0} a_{1} a_{2} b_{0} - 
   6 a_{2}^{2} b_{0}^{2} - a_{0} a_{1}^{2} b_{1} + 
   4 a_{0}^{2} a_{2} b_{1} + 4 a_{1} a_{2} b_{0} b_{1} \\
   &\qquad{}-2 a_{0} a_{2} b_{1}^{2} + a_{2} b_{1}^{3} - 
   2 a_{1}^{2} b_{0} b_{2} + 4 a_{0} a_{2} b_{0} b_{2} - 
   4 a_{2} b_{0} b_{1} b_{2} - a_{1} b_{1}^{2} b_{2} \\ 
   &\qquad{}+2 a_{0}^{2} b_{2}^{2} + 4 a_{1} b_{0} b_{2}^{2}, \\
   \rho(a,b)\s_2(\f) &=- a_{0}^{2} a_{1}^{2}  + 4 a_{0}^{3} a_{2} - 
  2 a_{1}^{3} b_{0} + 10 a_{0} a_{1} a_{2} b_{0} + 
  12 a_{2}^{2} b_{0}^{2} - 4 a_{0}^{2} a_{2} b_{1} \\ 
  &\qquad{}-7 a_{1} a_{2} b_{0} b_{1} - a_{1}^{2} b_{1}^{2} + 
  5 a_{0} a_{2} b_{1}^{2} - 2 a_{2} b_{1}^{3} + 
  2 a_{0}^{2} a_{1} b_{2} + 5 a_{1}^{2} b_{0} b_{2} \\
  &\qquad{}-4 a_{0} a_{2} b_{0} b_{2} - a_{0} a_{1} b_{1} b_{2} + 
  10 a_{2} b_{0} b_{1} b_{2} - 4 a_{1} b_{0} b_{2}^{2} + 
  2 a_{0} b_{1} b_{2}^{2} - b_{1}^{2} b_{2}^{2} \\
  &\qquad{}+4 b_{0} b_{2}^{3}. \\
  \endalign
$$
These formulas make the map $\Rat_2\to\M_2$ completely explicit using
the identifications $\Rat_2\subset\PP^5$ and
$\M_2\cong\Spec\ZZ[\s_1,\s_2]$.
\endremark

We will prove Theorem~\5.1 in a number of steps. One of the tools we
will use is a set of normal forms for rational maps of degree two
modulo $\SL_2$-conjugation. Normal forms are typically created by
moving fixed, periodic, and/or critical points into specified
locations, see for example [\MILN, appendix~C]. We will take the same
approach, but some care is needed because ultimately we will be
working over rings and fields which may have finite characteristic,
including characteristic~2. So for example, we will not want to use,
either implicitly or explicitly, the ``fact'' that a rational map of
degree~2 has exactly two critical points, since in characteristic~2
a map of degree two either has one critical point, or else it is
inseparable and every point is critical.

\remark{Remark}
The relation $\s_1=\s_3+2$, that is $\m_1+\m_2+\m_3=\m_1\m_2\m_3+2$,
implies the formal identities
$$
  (\m_1-1)^2=(\m_1\m_2-1)(\m_1\m_3-1)\qquad\hbox{and}\qquad
  (\m_2-1)^2=(\m_2\m_1-1)(\m_2\m_3-1).
$$
In particular, over any field (or even over any reduced ring), the
condition $\m_1\m_2=1$ is equivalent to $\m_1=\m_2=1$.
\endremark

\proclaim{Normal Forms Lemma \5.3}
Let $\O$ be an algebraically closed field of characteristic $p$, and
let $\xi\in\M_2(\O)$ have multipliers
$\m_1,\m_2,\m_3$. For any $\f(z)=F(z)/G(z)$ with $F,G\in\O[z]$,
we will write $\r$ for the resultant $\Res(F,G)$.
\Part{(i)}
If $\m_1\m_2\ne1$, then
$$
  \f(z)={z^2+\m_1z\over\m_2z+1}\in\xi
$$
with $\r=1-\m_1\m_2$. 
\Part{(ii)}
If $\m_1\ne0$, then there is a $\b\in\O$ satisfying
$$
  \b^2=\left(1-{2\over\m_1}\right)^2-\m_2\m_3
$$
such that
$$
  \f(z)={1\over\m_1}\left(z+{1\over z}\right)+\b\in\xi
$$
and $\r=\m_1^2$. 
\endproclaim
\demo{Proof}
We start with any rational map
$$
  \f(z)={a_0z^2+a_1z+a_2\over b_0z^2+b_1z+b_2}\in\xi
$$
and make coordinate changes to put~$\f$ into the desired form.
\Part{(i)} 
As noted above, the condition $\m_1\m_2\ne1$ implies that
$\m_1\ne1$ and $\m_2\ne1$. The associated fixed points thus have
multiplicity~1, so they must be distinct. Making a change of
variables, we can move them to~$0$ and~$\infty$ respectively. This
means that $\f(z)=(a_0z^2+a_1z)/(b_1z+b_2)\in\xi$ with $a_0b_2\ne0$,
so we can dehomogenize by setting $a_0=1$. Taking derivatives, we find
that
$$
  \f'(0)={a_1\over b_2}=\m_1 \qquad\hbox{and}\qquad
  \f'(\infty)=b_1=\m_2,
$$
so $\f$ has the form $\f(z)=(z^2+b_2\m_1z)/(\m_2z+b2)$. Finally,
$b_2^{-1}\f(b_2z)$ puts~$\f$ into the desired form, and one easily
verifies that the resultant is $1-\m_1\m_2$.
\Part{(ii)}
The assumption that $\m_1\ne0$ means that the associated fixed point
is not critical. We move this fixed point to~$\infty$, which
forces $b_0=0$, and then $a_0\ne0$, so we can dehomogenize by setting
$a_0=1$. Further, $\f'(\infty)=b_1=\m_1$. Next we observe that
$\f^{-1}(\infty)$ consists of~$\infty$ and one other point. This
follows from the fact that the multiplier at~$\infty$ is non-zero, or
we can just note that $\f(-b_2/\m_1)=\infty$. In any case, we use the
change of variables $z\mapsto z-b_2/\m_1$ to move this point to~$0$,
which puts~$\f$ in the form
$$
  \f(z)={z^2+a_1z+a_2\over\m_1z}.
$$
Note that $a_2\ne0$, so the final variable change
$z\mapsto\sqrt{a_2}z$ puts~$\f$ into the desired form with
$\b=a_1/\m_1\sqrt{a_2}\,$. The resultant is easily computed to equal
$\r=\m_1^2$, and with a bit more effort one computes the multiplier
$$
  \s_3=\m_1\m_2\m_3=\m_2-\m_2\b^2-4+4/\m_1.
$$
Solving for $\b^2$ completes the proof of the lemma.
\enddemo

Using these normal forms, it is not hard to show that the map
$\M_2\to\AA^2$ is bijective on geometric points. This may be compared
with [\MILN, lemma~3.1], where the same result is proven over~$\CC$
in essentially the same way.

\proclaim{Lemma \5.4}
Let $\Omega$ be an algebraically closed field. Then the map
$$
  (\s_1,\s_2):\M_2(\Omega)\longrightarrow \O^2
$$
is a bijection (of sets).
\endproclaim
\demo{Proof}
Let $\xi,\xi'\in\M_2(\O)$ have the same image in~$\O^2$.
The set of multipliers is determined by the values of~$\s_1$
and~$\s_2$, since the multipliers (with multiplicity) are the roots
of the polynomial
$$ 
  T^3-\s_1T^2+\s_2T-(\s_1-2).
$$
Hence~$\xi,\xi'$ have the same multiplier systems, say
$\{\m_1,\m_2,\m_3\}$. We consider two cases.
\par
First, suppose that $\m_1\m_2\ne1$. Then the Normal
Forms Lemma~\5.3(i) tells us that the rational map
$$
  \f(z)={z^2+\m_1z\over\m_2z+1}
$$
is in both~$\xi$ and~$\xi'$, so $\xi=\xi'$. 
\par
Second, suppose that
$\m_1\m_2=1$. The Normal Forms Lemma~\5.3(ii) says that~$\xi$
and~$\xi'$ each contain  a map
$$
  \f(z)={1\over\m_1}\left(z+{1\over z}\right)+\b\in\xi
$$
for some $\b\in\O$ satisfying
$$
  \b^2=\left(1-{2\over\m_1}\right)^2-\m_2\m_3.
$$
However, by an earlier remark, the condition $\m_1\m_2=1$ actually
implies that $\m_1=\m_2=1$, and then $\s_1=\s_3+2$ shows that also
$\m_3=1$. Hence $\b^2=0$, so $\b=0$. Thus~$\xi$
and~$\xi'$ both contain~$z+1/z$, so $\xi=\xi'$. This completes the
verification that the map $\M_2(\O)\to\O^2$ is injective.
\par
To see that the map is surjective, we take any
$(\a_1,\a_2)\in\Omega^2$, and we let
$\m_1,\m_2,\m_3\in\Omega$ be the three roots (with multiplicity) of
the polynomial
$$
  T^3-\a_1T^2+\a_2T-(\a_1-2).
$$
If any $\m_1\m_2\ne1$, then the rational map
$$
  \f(z)={z^2+\m_1z\over\m_2z+1}
$$
has degree two, multipliers~$\m_1,\m_2,\m_3$, and hence invariants
$\s_1(\f)=\a_1$ and $\s_2(\f)=\a_2$. Similarly, if $\m_1\m_2=1$,
then it follows as usual that $\m_1=\m_2=\m_3=1$, so we need merely
observe that the map $z+1/z$ has multipliers $\{1,1,1\}$ and
invariants $\s_1=3=\a_1$ and $\s_2=3=\a_2$. This proves that the map
$\M_2(\O)\to\O^2$ is surjective, which completes the proof of the
lemma.
\enddemo

Before proceeding further, we want to note that the mere fact that
$\M_2\to\AA^2$ is bijective on geometric points (i.e.,
$\M_2(\O)=\AA^2(\O)$ for algebraically closed fields~$\O$) does not
imply that the map~$\M_2\to\AA^2$ is an isomorphism.  The are two
possible problems. First, if~$\O$ has positive characteristic, then
an inseparable map may be bijective on points, yet not be an
isomorphism. Second, even in characteristic~$0$, there are morphisms
which are bijective on geometric points, yet have no inverse. A simple
example is the map of $\AA^1$ onto the the cuspidal cubic
$y^2=x^3$ via the map $t\mapsto(t^2,t^3)$. 
The next step in the proof of Theorem~\5.1 will be to show that the
map $\M_2\to\AA^2$ is proper. 

\proclaim{Lemma \5.5}
The map $\M_2\to\Spec\ZZ[\s_1,\s_2]=\AA^2_\ZZ$ is a proper morphism.
\endproclaim
\demo{Proof}
Let $F:\M^2\to\AA^2_\ZZ$ be the given map. We know from general
principles that~$\M_2$ and~$\AA^2_\ZZ$ are separable over~$\ZZ$,
so~$F$ is separable. Further,~$\M_2$ is of finite type of~$\ZZ$,
so~$\M_2$ is Noetherian and~$F$ is of finite type. Hence we may use
the valuative criterion \cite{\HART,~II.4.7} to check that~$F$ is
proper.
\par
Let
$$\alignedat2
  R &=\text{a discrete valuation ring,}\qquad & T &=\Spec(R), \\
  K &=\text{the fraction field of $R$,} & U &=\Spec(K),\\
  \endalignedat
$$
and suppose we are given a commutative diagram
$$\CD
  U=\Spec(K) @>>> \M_2 \\
  @VViV           @VVFV \\
  T=\Spec(R) @>>> \AA^2_\ZZ\rlap{.} \\
  \endCD
  \tag \neweqno{valuative square}
$$
We need to find a map $T\to\M_2$ making the diagram commute.
\par
Let $\Mbar_2=\M_2^{ss}$ be the proper scheme containing~$\M_2$ as
described in Theorem~\2.1. Since we are working with maps of degree~2,
Corollary~\2.3 says that $\M_2^{ss}=\M_2^{s}$, so
$\Mbar_2$ is actually the $\SL_2$-quotient of the stable points
in~$\PP^5$, but for our purposes it suffices to
know that there is a certain $\SL_2$-stable subset $(\PP^5)^{ss}$
of $\PP^5$ which contains $\Rat_d$ and whose $\SL_2$-quotient
$\Mbar_2$ exists and is proper over~$\ZZ$. In other words, we have a
commutative diagram
$$\CD
        @.  \Rat_2 @>>> (\PP^5)^{ss} \\
  @.          @VVV        @VVV   \\
  U    @>>>  \M_2  @>>>  \Mbar_2 \\
  @VVV       @VVV         @VVV   \\
  T    @>>>  \AA^2_\ZZ @>>> \Spec(\ZZ). \\
  \endCD
$$
The map $\Mbar_2\to\Spec(\ZZ)$ is proper, so the valuative criterion
implies that there is a map $T\to\Mbar_2$ making the diagram commute.
So we just need to show that the image of this map lies in~$\M_2$,
since this will give a map $T\to\M_2$, and then the separability
of~$\AA^2$ over~$\ZZ$ will imply that $T\to\M_2$ commutes with the
maps in the left-hand square.
\par
(To verify this last assertion, we label some of the maps in the above
diagram as
$$\CD
  U    @>\a>>  \M_2   \\
  @VViV       @VV{\s}V  \\
  T    @>\b>>  \AA^2_\ZZ @>\p>> \Spec(\ZZ). \\
  \endCD
$$
Now suppose that we have constructed a map $\g:T\to\M_2$ with
$\g\circ i=\a$, but that we only know that
$\p\circ\s\circ\g=\p\circ\b$. We want to show that
$\s\circ\g=\b$. Of course, we do know that $\s\circ\a=\b\circ i$.
Consider the commutative square
$$\CD
  U @>\s\circ\g\circ i=\s\circ\a=\b\circ i>> \AA^2 \\
  @VViV  @VV{\p}V \\
  T @>\p\circ\s\circ\g=\p\circ\b>> \Spec(\ZZ). \\
  \endCD
$$
Notice that both of the maps $\s\circ\g:T\to\AA^2$ and $\b:T\to\AA^2$
commute with this square. Hence the separability of
$\p:\AA^2\to\Spec(\ZZ)$ implies the desired equality $\s\circ\g=\b$.)
\par
We observe that we are free to replace~$K$ by a finite extension~$K'$
and~$R$ with its integral closure~$R'$ in~$K'$. This is true because
if we can prove that the  map $T'=\Spec(R')\to\Mbar_2$ has image
in~$\M_2$, then the same will be true for $T\to\Mbar_2$, since
$T'\to T$ is surjective. 
\par
The given map $\xi:U\to\M_2$ is  a $K$-valued point
$\xi\in\M_2(K)$. 
This $\SL_2$-equivalence class of rational maps
has invariants $\s_1,\s_2,\s_3\in K$ and multipliers
$\m_1,\m_2,\m_3\in\Kbar$ as usual. Replacing~$K$ by a finite
extension, we will assume that $\m_1,\m_2,\m_3\in K$.
Further, the commutativity of the 
diagram~(\referto{valuative square}) tells us that~$\s_1$ and~$\s_2$
are in~$R$. Hence $\m_1,\m_2,\m_3$ are also in~$R$, since they are
roots of the monic polynomial with coefficients in~$R$,
$$
  T^3-\s_1T^2+\s_2T-(\s_1-2),
$$
and~$R$ is integrally closed.
Let $\gM$ be the maximal ideal in the valuation ring~$R$. We now
consider two cases.
\par
First, suppose that $\m_1\m_2\not\equiv1\pmod{\gM}$. Then certainly
$\m_1\m_2\ne1$, so the Normal Forms Lemma~\5.3(i) tells us that (after
another finite extension of~$K$) we can find a map
$$
  \f(z)={z^2+\m_1z\over\m_2z+1}
$$
in the equivalence class of maps~$\xi$. In other words,~$\f$ is a
point in~$\Rat_2(K)$ lifting~$\xi$. Recall that~$\Rat_2$ is the
affine open subset of~$\PP^5$ given by
$$
  \bigl\{[a_0,a_1,a_2,b_0,b_1,b_2]\in\PP^5\,:\,
  \Res(a_0X^2+a_1XY+a_2Y^2,b_0X^2+b_1XY+b_2Y^2)\ne0\bigr\}.
$$
Thus a point $\psi\in\Rat_2(K)\hookrightarrow\PP^5(R)$ will lie
in~$\Rat_2(R)$ if and only if it has the form
$$\gather
  \psi(z)={a_0z^2+a_1z+a_2\over b_0z^2+b_1z+b_2}
    \quad\text{with $a_0,a_1,a_2,b_0,b_1,b_2\in R$ and} \\
  \Res(a_0X^2+a_1XY+a_2Y^2,b_0X^2+b_1XY+b_2Y^2)
      \not\equiv0\pmod{\gM}. \\
  \endgather
$$
The map~$\f(z)$ listed above corresponds to the point
$[1,\m_1,0,0,\m_2,1]\in\PP^5(R)$, and its resultant is $\m_1\m_2-1$.
By assumption, $\m_1\m_2-1\not\equiv0\pmod{\gM}$, so~$\f$ lies
in~$\Rat_2(R)$, and hence~$\xi$ lies in~$\M_2(R)$.
\par
For the second case, we suppose that $\m_1\m_2\equiv1\pmod{\gM}$.
In particular, $\m_1\ne0$, so the Normal Forms Lemma~\5.3(ii) says
that (after a finite extension of~$K$) there is a map
$$
  \f(z)={1\over\m_1}\left(z+{1\over z}\right)+\b\in\xi,
$$
where~$\b$ satisfies
$$
  \b^2=\left(1-{2\over\m_1}\right)^2-\m_2\m_3.
$$
Again extending~$K$, we may assume that~$\b\in K$. Further, the
assumption that $\m_1\m_2\equiv1\pmod{\gM}$ means that~$\m_1$ is a
unit (i.e.,~$\m_1\in R^*$), so we see that~$\b\in R$. As above, this
map~$\f$ corresponds to the point $[1,\b\m_1,1,0,\m_1,1]\in\PP^2(R)$
having resultant~$\m_1^2$. We know that~$\m_1$ is a unit, so
$\m_1^2\not\equiv0\pmod{\gM}$. Hence~$\f$ lies in~$\Rat_2(R)$, which
proves that~$\xi$ lies in~$\M_2(R)$.
\par
This completes the proof that there is a map $T=\Spec(R)\to\M_2$
making the diagram~(\referto{valuative square}) commute. By the
valuative criterion for properness \cite{\HART,~II.4.7}, we conclude
that the map $(\s_1,\s_2):\M_2\to\AA^2_\ZZ$ is proper, which
completes the proof of Lemma~\5.5.
\enddemo

We now know that the map $\M_2\to\AA^2_\ZZ$ is proper. However,
both~$\M_2$ and~$\AA^2_\ZZ$ are affine varieties. Intuitively, the
(geometric) fibers of the map are both affine and complete, which
should imply  they they consist of a finite set of points. The
following generalization of \cite{\HART,~exer.~II.4.6} makes this
intuition precise. It will be used to show that the
map~$\M_2\to\AA^2_\ZZ$ is finite.

\proclaim{Lemma \5.6}
Let $X$ and $Y$ be affine integral schemes, and let $F:X\to Y$ be a
dominant proper morphism of finite type. Then~$F$ is a finite
morphism.
\endproclaim
\demo{Proof}
Let $X=\Spec(A)$ and $Y=\Spec(B)$. Since~$X$ and~$Y$ are integral
schemes,~$A$ and~$B$ are integral domains. We let~$K_A$ and~$K_B$ be
the fraction fields of~$A$ and~$B$ respectively. Then
$F:X\to Y$ induces a homomorphism $B\to A$, and the fact that~$F$ is
dominant means that this homomorphism is injective. So we also get an
injection $K_B\to K_A$.
\par
Now let~$B'\subset K_A$ be any valuation ring of~$K_A$ containing the
image of~$B$. By definition of valuation ring, every $x\in K_A$
satisfies either $x\in B'$ or $x^{-1}\in B'$. (See
\cite{\LANG,~XII~\S4}.) In particular,~$K_A$ is the fraction field
of~$B'$. Now consider the commutative diagrams
$$\CD
    A @>>> K_A \\  @AAA  @AAA \\ B @>>> B' \\
  \endCD
  \qquad\qquad
  \CD
    X=\Spec(A) @<<< \Spec(K_A) \\
    @VVFV @VVV\\
    \Spec(B) @<<< \Spec(B')\rlap{.} \\
  \endCD
$$
We are given that~$F$ is proper, so the valuative criterion of
properness \cite{\HART,~II.4.7} tells us that there is a unique map
$\Spec(B')\to\Spec(A)$ making the right-hand diagram commute.
Equivalently, there is a unique homomorphism $A\to B'$ making the
left-hand diagram commute. This proves that every valuation ring
of~$K_A$ containing the image of~$B$ will also contain~$A$. It follows
from \cite{\LANG,~XII~\S4, prop.~4.9} or \cite{\HART,~II.4.11A}
that~$A$ is integral over~$B$. (That is, every element of~$A$
is the root of a monic polynomial in~$B[T]$.) But we are also given
that~$F$ is of finite type, which means that~$A$ is of finite type
over~$B$. Thus $A=B[a_1,\ldots,a_r]$ with each~$a_i$ integral
over~$B$, so~$A$ is a finitely generated $B$-module. Therefore~$F$ is
finite.
\enddemo

Combining Lemmas~\5.4,~\5.5, and~\5.6 shows that the map
$\M_2\to\AA^2_\ZZ$ is a finite map which is bijective on geometric
points. One would expect that this should imply that the map is an
isomorphism, but there is still some work to do. In
characteristic~$p$, the Frobenius map is finite and bijective on
geometric points, yet is not an isomorphism; and the same is true of
the map $t\to(t^2,t^3)$ of~$\AA^1$ onto the twisted cubic. We will
need to use the fact that $\M_2\to\AA^2_\ZZ$ is a morphism of schemes
over~$\ZZ$ and the fact that the image~$\AA^2_\ZZ$ is non-singular.
The following general lemma is what we will need to complete the
proof of Theorem~\5.1.

\proclaim{Lemma \5.7}
Let $F:X\to Y$ be a morphism of schemes over~$\ZZ$. Suppose that the
following four conditions are true.
\Part{(i)}
$X$ is an integral scheme.
\Part{(ii)}
$Y$ is an integral normal scheme which is dominant over~$\ZZ$.
\Part{(iii)}
$F$ is a finite morphism.
\Part{(iv)}
$F$ induces a bijection on geometric points.
\par\noindent
Then $F$ is an isomorphism.
\endproclaim
\demo{Proof}
Before beginning the proof, we remind the reader what the
conditions~(i)--(iv) really mean. A scheme~$X$ is integral if and only
if it is reduced and irreducible~\cite{\HART,~II.3.1}. This implies
that the local ring of the generic point is a field, 
equal to the fraction field of~$A$ for any affine open subset
$\Spec(A)\subset X$~\cite{\HART,~II.3.6}. The fact that~$Y$ is normal
means that its local rings are integrally closed
domains~\cite{\HART,~exer.~II.3.8}. The map~$F$ induces a natural map 
on $S$-valued points, $X(S)\to Y(S)$, for any scheme~$S$. Thus if
$f:S\to X$ is in~$X(S)$, then $F\circ f\in Y(S)$. Condition~(iv) says
that this map is a bijection whenever $S=\Spec(\O)$ for an
algebraically closed field~$\O$.
\par
We now begin the proof of Lemma~\5.7. Geometric points are dense
in~$Y$, so~(iv) implies that~$F$ is dominant (i.e.,~$F(X)$ is dense
in~$Y$). Hence~$F$ induces a map of function fields
$F^*:K(Y)\hookrightarrow K(X)$. (In fact, since~$F$ is finite
from~(iii), it follows that~$F$ is
closed~\cite{\HART,~ex.~II.3.5(b)}, so~$F(X)$ is dense and closed,
so~$F(X)=Y$. But we won't need to know this stronger fact.)
\par
Take any affine open subset $\Spec(B)\subset Y$. Using the fact~(iii)
that~$F$ is finite, we find that 
$F^{-1}\bigl(\Spec(B)\bigr)=\Spec(A)$, where~$A$ is a finitely
generated $B$-module~\cite{\HART,~exer.~II.3.4}. Replacing~$X$
and~$Y$ by $\Spec(A)$ and $\Spec(B)$, we may assume that~$X$ and~$Y$
are affine. Note that~$A$ and~$B$ are integral domains, since~$X$
and~$Y$ are integral schemes from~(i) and~(ii). Let~$K_A$ and~$K_B$
be the fraction fields of~$A$ and~$B$ respectively. We have
commutative diagrams
$$
  \CD
    X=\Spec(A) @<<< \Spec(K_A) \\
    @VVFV   @VVFV \\
    Y=\Spec(B) @<<< \Spec(K_B) \\
  \endCD
  \qquad\qquad
  \CD
    A @>>> K_A \\
    @AAF^*A  @AAF^*A \\
    B @>>> K_B\rlap{.} \\
  \endCD
$$
\par
Notice that if~$\O$ is any (algebraically closed) field, then
$$
  X(\O)=\Mor(\Spec\O,X)=\Mor(\Spec\O,\Spec A)
  =\Hom(A,\O)=\Hom(K_A,\O),
$$
where the last equality is true becasue~$K_A$ is the fraction field
of the integral domain~$A$, and similarly $Y(\O)=\Hom(K_B,\O)$. Then
the map  $X(\O)\to Y(\O)$ induced by~$F$ is given by
$$
  \Hom(K_A,\O)\longrightarrow\Hom(K_B,\O),\qquad
  f\longmapsto f\circ F^*.
$$
Now it is a standard fact from the theory of fields that if~$\O$ is
algebraically closed, then any $g\in\Hom(K_B,\O)$ can be lifted to an
element of $\Hom(K_A,\O)$ in exactly $[K_A:F^*K_B]_s$ ways, where the
subscript~$s$ denotes the separable degree. (See
\cite{\LANG,~VII~\S4}.) Condition~(iv) tells us that $X(\O)\to Y(\O)$
is bijective, so we conclude that $[K_A:F^*K_B]_s=1$. However, the
assumption~(ii) that~$Y$ is dominant over~$\Spec(\ZZ)$ implies
that~$B$, and hence also~$K_B$, have characteristic~$0$. 
So the separable degree is the actual degree, $[K_A:F^*K_B]=1$, and
hence $F^*:K_B\to K_A$ is an isomorphism.
\par
We now have the commutative diagram
$$\CD
  A @>>> K_A \\
  @AAF^*A  @| \\
  B @>>> K_B\rlap{.} \\
  \endCD
$$
Further,~$A$ is integral over~$B$ from~(iii), and~$B$ is integrally
closed in~$K_B$ from~(ii). But~$K_A=F^*(K_B)$, so~$F^*(B)$ is
integrally closed in~$K_A$. This gives the inclusions
$$
  B @>F^*>> A 
  \subset \text{(integral closure of $B$ in $K_A$)} 
  = F^*(B).
$$
Hence $F^*:B\to A$ is an isomorphism, which completes the proof
that~$F$ is an isomorphism.
\enddemo

We now have all of the pieces to complete the proof of Theorem~\5.1.

\demo{Proof of Theorem \5.1}
We will denote by $\s:\M_2\to\AA^2_\ZZ$ the morphism induced by the
inclusion of $\ZZ[\s_1,\s_2]$ into the affine coordinate ring
of~$\M_2$. Both~$\M_2$ and~$\AA^2_\ZZ$ are affine integral schemes,
the former from Theorem~\2.1, and the latter trivially. The
map~$\s$ is bijective on geometric points from Lemma~\5.4, so it
must be dominant. The map~$\s$ is of finite type, since it is a
$\ZZ$-morphism, and~$\M_2$ is actually of finite type over~$\ZZ$ from
Theorem~\2.1. Finally,~$\s$ is a proper morphism from Lemma~\5.5.
It follows from Lemma~\5.6 that~$\s$ is a finite morphism.
\par
We want to apply Lemma~\5.7 to~$\s$, so we have to check the four
conditions in Lemma~\5.7. First,~$\M_2$ is an integral scheme from
Theorem~\2.1. Second, it is easy to see that~$\AA^2_\ZZ$ is an
integral normal scheme and is dominant over~$\ZZ$. Third,~$\s$ is a
finite morphism from the previous paragraph. Fourth,~$\s$ induces a
bijection on geometric points from Lemma~\5.4. Hence we can apply
Lemma~\5.7 to conclude that~$\s$ is an isomorphism.
\enddemo

\demo{Proof of Corollary \5.2}
In general, the affine coordinate ring of $\M_d$ is the ring
of $\SL_2$-invariant functions on $\Rat_d$ (see
Theorem~\2.1), while the affine coordinate ring of
$\AA^2_\ZZ=\Spec\ZZ[\s_1,\s_2]$ is precisely $\ZZ[\s_1,\s_2]$. Now
the isomorphism $\M_2\cong\AA^2_\ZZ$ from Theorem~\5.1 shows
that~$\M_2$ and~$\AA^2_\ZZ$ have the same affine coordinate
rings, which gives the first part of the corollary. The second part
is immediate, since   the higher order invariants
$\s_i^{(n)}$ are in the affine coordinate ring of~$\M_2 $.
\enddemo


\head \S\6. The completion $\M_2^s$ of $\M_2$
\endhead

In this section we will prove that the stable completion $\M_2^s$
of~$\M_2$ has a very simple structure as described in the following
theorem.

\proclaim{Theorem \6.1}
The isomorphism $(\s_1,\s_2):\M_2\cong\AA^2$ extends to an
isomorphism $\s:\M_2^s\cong\PP^2$ of schemes over~$\ZZ$.
\endproclaim

\remark{Remark}
Milnor \cite{\MILN, section~4} uses the identification
$(\s_1,\s_2):\M_2(\CC)\cong\CC^2$ to study the completion $\PP^2(\CC)$
of $\CC^2$. He shows that the extra points at infinity in~$\PP^2(\CC)$
correspond to linear and constant maps which can be thought of as
degenerate quadratic maps. This provides a natural completion
of~$\M_2(\CC)$ which is isomorphic to~$\PP^2(\CC)$, but unfortunately
it does not immediately imply that our completion~$\M_2^s(\CC)$ is
isomorphic to~$\PP^2(\CC)$. The difficulty is that Milnor implicitly
defines a degenerating family of maps $\f_t\in\Rat_d$ to be a family
for which (at least) one of~$\s_1(\f_t)$ or~$\s_2(\f_t)$ tends to
infinity as $t\to t_0$, but there is no a priori reason that the
$\s_i(\f_t)$'s might not approach some indeterminate
form~$\frac{0}{0}$.  For example, the family of maps
$$\gather
  \f_{(t_1,t_2)}(z)=\frac{t_1z^2+2z}{t_2z+1} \\
  \noalign{\flushpar\text{over the $(t_1,t_2)$ plane satisfies}}
  \s_1(\f_{(t_1,t_2)})=\frac{2(t_1^2-2t_1t_2-t_2^2)}{t_1(t_1-2t_2)}
  \quad\text{and}\quad
  \s_2(\f_{(t_1,t_2)})=\frac{-5t_2^2}{t_1(t_1-2t_2)},\\
  \endgather
$$
so neither of the limits
$$
  \lim_{(t_1,t_2)\to(0,0)} \s_1(\f_{(t_1,t_2)})
  \qquad\text{or}\qquad
  \lim_{(t_1,t_2)\to(0,0)} \s_2(\f_{(t_1,t_2)})
$$
exists. Milnor's completion corresponds to adding maps which 
satisfy what one might call a $(\s_1,\s_2)$-stability condition. During
the proof of Theorem~\6.1, we will verify that $(\s_1,\s_2)$-stability
is the same as the stability criterion from geometric invariant
theory used to define the stable sets~$(\PP^5)^s$ and~$\M_2^s$.
\endremark

\proclaim{Lemma \6.2}
Let $\O$ be an algebraically closed field, and let
$\f\in\bigl((\PP^5)^s\setminus\Rat_2\bigr)(\O)$. That
is,~$\f=[F_a,F_b]$ is in the stable locus, but~$\f$ is not in~$\Rat_2$
because the resultant $\Res(F_a,F_b)$ vanishes. Then there exists an
$f\in\SL_2(\O)$ so that
$$
  \f^f(X,Y) = [AXY,XY+BY^2]
  \qquad\text{for some $[A,B]\in\PP^1(\O)$.}
$$
Further, the homogeneous pair~$[A,B]$ is uniquely determined by the
conjugacy class $\<\f\>$ up to reversing the roles of~$A$ and~$B$.
\par
In other words, there is a well-defined bijection
$$
  \PP^1(\O)/\Sym_2 \longrightarrow
        (\M_2^s\setminus\M_2)(\O),
   \quad\text{induced by}\quad
   [A,B] \longmapsto [AXY,XY+BY^2],
$$
where the symmetric group on two letters~$\Sym_2$ acts on~$\PP^1$ by
interchanging the coordinates.
\endproclaim
\demo{Proof}
The assumption that~$\f=[F_a,F_b]$ is not in~$\Rat_2(\O)$ means
that~$F_a$ and~$F_b$ have a common root in~$\PP^1(\O)$. Making an
appropriate conjugation, we may move the common root to~$[1,0]$,
so~$\f$ looks like
$$
  \f=[a_1XY+a_2Y^2,b_1XY+b_2Y^2].
$$
According to Proposition~\2.2, the stability of~$\f$ 
(i.e.,~$\f\in(\PP^5)^s$) implies that $b_1\ne0$.
\par
Of course, we are still free to conjugate by elements of $\SL_2(\O)$
which fix~$[1,0]$. First we will conjugate by the matrix
$f=\left({1\atop0}\,{\b\atop1}\right)$. This gives
$$
  \f^f=\bigl[ (a_1-b_1\b )XY - (b_1\b^2+b_2\b-a_1\b-a_2)Y^2,
     b_1XY+(b_2+b_1\b)Y^2 \bigr].
$$
We know that~$b_1\ne0$, so we can take~$\b$ to be either of the roots
of
$$ 
  b_1\b^2+b_2\b-a_1\b-a_2=0
$$
and dehomogenize by setting $b_1=1$ to obtain
$$
  \f^f=[AXY,XY+BY^2]\qquad\text{with $A,B\in\O$.}
$$
If either~$A$ or~$B$ is non-zero, this is the desired form. But if
$A=B=0$, so $\f^f=[0,XY]$, then conjugation by
$\left({0\atop1}\,{-1\atop\phantom{-}0}\right)$ would give the form
$[XY,0]$, and from Proposition~\2.2 this would contradict the
stability of~$\f$. This shows that after conjugation, we can always
put~$\f$ into the desired form.
\par
It remains to determine to what extent  the form $[AXY,XY+BY^2]$ is
unique.  Some algebra and a  case-by-case analysis shows that
conjugation by the matrix $\left({\a\atop\g}\,{\b\atop\d}\right)$
preserves this form in exactly the following cases, where~$u$ denotes
an arbitrary element $u\in\O^*$:
$$\xalignat{3}
  f&=\pmatrix u^{-1}&0\\ 0&u\\ \endpmatrix,
    &\f^f&=[u^2AXY,XY+u^2BXY]
    &&\text{(any $A,B$)} \\
  f&=\pmatrix u^{-1}&u(A-B)\\ 0&u\\ \endpmatrix,
    &\f^f&=[u^2BXY,XY+u^2AXY]
    &&\text{(any $A,B$)} \\
  f&=\pmatrix u^{-1}&uB\\ -(uB)^{-1}&0\\ \endpmatrix,
    &\f^f&=[-u^2BXY,XY]
    &&\text{(if $A=0$, $B\ne0$)} \\
  f&=\pmatrix u^{-1}&0\\ -(uB)^{-1}&u\\ \endpmatrix,
    &\f^f&=[0,XY-u^2BY^2]
    &&\text{(if $A=0$, $B\ne0$)} \\
  f&=\pmatrix 0&uA\\ -(uA)^{-1}&0\\ \endpmatrix,
    &\f^f&=[-u^2AXY,XY]
    &&\text{(if $B=0$, $A\ne0$)} \\
  f&=\pmatrix 0&uA\\ -(uA)^{-1}&u\\ \endpmatrix,
    &\f^f&=[0,XY-u^2AY^2]
    &&\text{(if $B=0$, $A\ne0$)} \\
  \endxalignat
$$
It follows that two forms $[AXY,XY+BY^2]$ and $[A'XY,XY+B'Y^2]$ are
$\SL_2(\O)$-equivalent if and only if there is a $\l\in\O$ such that
either
$$
  (A',B') = (\l A,\l B) \qquad\text{or}\qquad (A',B') = (\l B,\l A).
$$
So the set of forms, up to $\SL_2(\O)$-equivalence, is in one-to-one
correspondence with the quotient space $\PP^1(\O)/\Sym_2$.
\enddemo

\proclaim{Lemma \6.3}
\rom{(a)}
Let $R$ be a discrete valuation ring with fraction field~$K$ and
residue field~$k$. Let $\psi:\Spec R\to\M_2^s$ be a morphism. Then the
point
$$
  \bigl[1,\s_1(\psi),\s_2(\psi)\bigr]\widetilde{\phantom{(((}}
  \in\PP^2(k)
$$
depends only on the image of the special fiber $\psi_k$. \rom(The
tilde indicates the natural reduction map $\PP^2(K)\to\PP^2(k)$.\rom)
\Part{(b)}
The map $(\s_1,\s_2):\M_2\to\AA^2$ induces a birational $\ZZ$-morphism
$$
  \s=[1,\s_1,\s_2]:\M_2^s\to\PP^2.
$$
\endproclaim
\demo{Proof}
(a)
Note that it is crucial that the special fiber~$\psi_k$ of the
family is assumed to be stable, since the example described in the
remark above shows that the result is false without the stability
assumption.
\par
Our first step will be to lift the map from $\M_2^s$ to $(\PP^5)^s$.
To do this, we first replace~$R$ by its strict Henselization. This
means that the residue field~$k$ is separably closed and that~$R$
satisfies Hensel's lemma. (See~\cite{\ATAEC,~IV~\S6} or
\cite{\BOSCH} for information about Henselizations.) 
\par
Next we observe that the map
$(\PP^5)^s\to\M_2^s$ is a smooth morphism. Intuitively, this is true
because it is a geometric quotient map whose fibers are isomorphic to
the smooth scheme~$\SL_2$. More precisely, we begin by
using~\cite{\ALTK, corollary~VII.1.9}. This says that it suffices to
check that the map over each point of $\Spec\ZZ$ is smooth. In other
words, we need to check that the maps
$$
  (\PP^5)^s\times\Spec\FF\longrightarrow\M_2^s\times\Spec\FF
$$
are smooth, where~$\FF$ is either $\ZZ/p\ZZ$ or~$\QQ$. This reduces the
problem to morphisms over a field. Next we apply
\cite{\HART,~III.10.2}, which says that it suffices to check that the
fiber over each point of~$\M_2^s(\bar\FF)$ is regular. But as noted
above, each such fiber is isomorphic to $\SL_2/\bar\FF$. This
completes the verification that the morphism $(\PP^5)^s\to\M_2^s$ is
smooth. Now the lifting property of Henselian rings \cite{\BOSCH} says
that the  map $(\PP^5)^s(R)\to\M_2^s(R)$ on
$R$-points is surjective, so we can lift~$\psi$. By abuse of notation,
we will also denote the lift by $\psi:\Spec R\to(\PP^5)^s$.
\par
If $\psi_k\in\Rat_2(k)$, then $\s_1(\psi)$ and $\s_2(\psi)$ are
already in~$R$, so the desired result follows from the trivial
computation
$$
  \bigl[1,\s_1(\psi),\s_2(\psi)\bigr]\widetilde{\phantom{(((}}
  =\bigl[1,\widetilde{\s_1(\psi)},\widetilde{\s_2(\psi)}\bigr]
  =\bigl[1,\s_1(\psi_k),\s_2(\psi_k)\bigr].
$$
\par
Suppose now that $\psi_k\notin\Rat_2(k)$. Of course, by assumption we
do know that $\psi_k\in(\PP^5)^s(k)$, so Lemma~\6.2 says that after
conjugation, we may assume that~$\psi_k$ has the form
$$
  \psi_k=[\tilde A_1 XY,XY+\tilde B_2 Y^2]
$$
for some $[\tilde A_1,\tilde B_2]\in\PP^1(k)$. Lifting~$\tilde A_1$
and~$\tilde B_2$ to elements $A_1,B_2\in R$, this means that we can
write~$\psi$ in the form
$$
  \psi=[A_0\p X^2+A_1 XY+A_2\p Y^2,B_0\p X^2+(1+B_1\p)XY+B_2 Y^2],
$$
where $A_0,A_1,A_2,B_0,B_1,B_2\in R$, at least one of~$A_1,B_2$ is
in~$R^*$, and $\p\in R$ is a uniformizer. Let
$$
  \r=\r(\psi)=\Res(A_0\p X^2+A_1 XY+A_2\p Y^2,
    B_0\p X^2+(1+B_1\p)XY+B_2 Y^2).
$$
Notice that $\tilde\r=0$. We claim that $\r(\psi)\s_1(\psi)$ and
$\r(\psi)\s_2(\psi)$ are both in~$R$ and that at least one of them is
in~$R^*$. To see this, we use the explicit formulas for~$\s_1$
and~$\s_2$ in section~\5. Substituting into these formulas and
reducing modulo~$\p$, we find that
$$ 
  \widetilde{\r(\psi)\s_1(\psi)} = -\widetilde{A_1B_2}
  \qquad\text{and}\qquad
  \widetilde{\r(\psi)\s_2(\psi)} = -\widetilde{A_1^2-B_2^2}.
$$
Hence
$$\align
  \bigl[1,\s_1(\psi),\s_2(\psi)\bigr]\widetilde{\phantom{(((}}
  &=\bigl[\r(\psi),\r(\psi)\s_1(\psi),
      \r(\psi)\s_2(\psi)\bigr]\widetilde{\phantom{(((}} \\
  &=\bigl[\widetilde{\r(\psi)},\widetilde{\r(\psi)\s_1(\psi)},
      \widetilde{\r(\psi)\s_2(\psi)}\bigr] \\
  &=[0,\widetilde{A_1B_2},\widetilde{A_1^2+B_2^2}]
  \in\PP^2(k). \\
  \endalign
$$
Notice that this is a well-defined point in~$\PP^2(k)$, since~$A_1$
and~$B_2$ are in~$R$ and at least one of them is a unit. Further, the
point clearly depends only on the special fiber $\psi_k=[\tilde A_1
XY,XY+\tilde B_2 Y^2]$. This completes the proof of Lemma~\6.3(a).
\Part{(b)}
The map $(\s_1,\s_2):\M_2\to\AA^2$ is an isomorphism from
Theorem~\5.1, so it certainly induces a birational map
$\M_2^s\to\PP^2$. We want to show that this map extends to a
morphism. This follows from~(a) and general principles. We briefly
sketch. If $F:X\to Y$ is a birational map with the property in~(a),
and if $x\in X$ is a closed point, we define~$F(x)$ as follows.
Take any discrete valuation ring~$R$ with fraction field~$K$
and residue field~$k(x)$ and any map $\psi:\Spec R\to X$ with
$\psi(\Spec k(x))=x$. Since~$F$ is birational, we can also assume that
$\psi(\Spec K)$ lies in the domain of~$F$. Then $F\circ\psi:\Spec
R\to Y$ extends to a morphism (assuming~$Y$ is regular), so we can
define $F(x)=(F\circ\psi)(\Spec k(x))$. The key here is that the
property described in~(a) says that the value of~$F(x)$ depends only
on~$x$, independent of~$\psi$, so~$F(x)$ is well-defined.
\enddemo

We now have the tools needed to complete the proof of Theorem~\6.1.

\demo{Proof of Theorem \6.1}
Lemma \6.3(b) says that there is a birational morphism
$$
  \s=[1,\s_1,\s_2]:\M_2^s\longrightarrow\PP^2.
$$
We claim that~$\s$ is a bijection on geometric points. Theorem~\5.1
tells us that~$\s$ is an isomorphism $\M_2\to\AA^2$, so we just need
to check the boundary. Let~$\O$ be an algebraically closed field, and
for any $A,B\in\O$, let $\f_{A,B}=[AXY,XY+BY^2]\in\PP^5(\O)$. Then
Lemma~\6.2 says that the map $[A,B]\to\f_{A,B}$ induces a bijection
$\PP^1(\O)/\Sym_2\to(\M_2^s\setminus\M_2)(\O)$, so we need to show
that
$$
  \PP^1(\O)/\Sym_2 \longrightarrow (\PP^2\setminus\AA^2)(\O),\qquad
  [A,B]\longmapsto \s(\f_{A,B}),
$$
is a bijection. Note that we cannot compute $\s(\f_{A,B})$ directly,
since $\s_1(\f_{A,B})$ and $\s_2(\f_{A,B})$ do not exist. However, if
we let~$\r$ denote the resultant form of two polynomials, then~$\r\s_1$
and~$\r\s_2$ will be defined at the point~$\f_{A,B}$. More precisely,
using the explicit formulas for~$\r\s_1$ and~$\r\s_2$ given in
section~\5, we find that
$$
  \r(\f_{A,B})=0,\qquad (\r\s_1)(\f_{A,B})=-AB,\qquad
  (\r\s_2)(\f_{A,B})=-A^2-B^2,
$$
so we are reduced to showing that the map
$$
  \PP^1(\O)/\Sym_2\longrightarrow (\PP^2\setminus\AA^2)(\O),\qquad
  [A,B]\longmapsto [0,AB,A^2+B^2]
$$
is bijective. This is an exercise which we leave for the reader.
\par
We now know that $\s:\M_2^s\to\PP^2$ is a birational morphism which
is bijective on geometric points. The fact that it is bijective on
geometric points certainly implies that it is quasi-finite (i.e., the
inverse image of any point is  a finite set of points).
Further,~$\M_2^s$ and~$\PP^2$ are both proper over~$\ZZ$, so it
follows from~\cite{\HART,~II.4.8(e)} that~$\s$ is a proper morphism.
Thus~$\s$ is quasi-finite and proper, so~\cite{\MILNE, chapter~I,
proposition~1.10} tells us that~$\s$ is finite. 
\par
To complete the proof of Theorem~\6.1, we merely need to observe that
we now know that $\s:\M_2^s\to\PP^2$ satisfies the four conditions in
Lemma~\5.7, and hence~$\s$ is an isomorphism.
(We remark that rather than using Lemma~\5.7, we could instead give a
direct proof that a finite birational morphism
$F:X\to Y$ of integral schemes with~$Y$ normal is an isomorphism. To
do this, we can replace~$X$ and~$Y$ by affines $\Spec A$ and $\Spec B$.
Then~$A$ is integral over~$B$, the fraction fields of~$A$ and~$B$
coincide, and~$B$ is integrally closed, so $A=B$.)
\enddemo


\remark{Acknowledgements}
I would like to thank Spencer Bloch, Ching-Li Chai, Joe Harris, Stephen
Lichtenbaum, and Curt McMullen for numerous helpful suggestions.
\endremark


\enddocument

\bye